\documentclass[a4paper, 11pt]{article}
\usepackage [a4paper,left=2.5cm,bottom=2.5cm,right=2.5cm,top=2.5cm]{geometry}
\usepackage[english]{babel}
\usepackage{amssymb}
\usepackage{amsmath,amsthm, dsfont}
\usepackage{subcaption}
\usepackage{comment}
\usepackage{tabularx,array}
\usepackage{graphicx, url}
\usepackage{enumitem} 
\usepackage{bm}
\usepackage{todonotes}
\usepackage{todonotes}
\usepackage{hyperref}
\hypersetup{colorlinks=true, linkcolor=blue,citecolor=gray}
\usepackage[nameinlink,capitalise]{cleveref}
\crefformat{equation}{#2(#1)#3}

\theoremstyle{plain}
\theoremstyle{plain}\newtheorem{assumption}{}

\newtheorem{theorem}{Theorem}[section]
\newtheorem{lemma}[theorem]{Lemma}

\newtheorem{proposition}[theorem]{Proposition}

\newtheorem{remark}[theorem]{Remark}
\newtheorem{definition}{Definition}

\newcommand{\be}{\begin{equation}}
\newcommand{\ee}{\end{equation}}

\newcommand{\bega}{\begin{equation}\begin{aligned}}
\newcommand{\eega}{\end{aligned}\end{equation}}

\newcommand{\E}{\mathbb{E}}

\newcommand{\R}{\mathbb{R}}

\renewcommand{\P}{\mathbb{P}}

\newcommand{\rev}{\color{black}}
\newcommand{\spa}{\color{black}}

\title{Kinetic interacting particle system: parameter estimation from complete and partial discrete observations}
\author{Chiara Amorino$^{*}$, Vytaut\.e Pilipauskait\.e$^{\dagger}$}

\date{\today}

\begin{document}
\maketitle

\footnote{
VP acknowledges support of the Independent Research Fund Denmark - Natural Sciences, grant DFF-10.46540/2032-00005B. \\

$^{*}$ Universitat Pompeu Fabra and Barcelona School of Economics, Barcelona, Spain

{$^\dagger$ Centrica Energy A/S, Aalborg, Denmark}

}

\begin{abstract}
In this paper, we study the estimation of drift and diffusion coefficients in a two-dimensional system of \(N\) interacting particles modeled by a degenerate stochastic differential equation. We consider both complete and partial discrete observation cases over a fixed time horizon \([0, T]\) and propose novel contrast functions for parameter estimation. In the partial observation scenario, we tackle the challenge posed by unobserved velocities by introducing a surrogate process based on the increments of the observed positions. This requires a modified contrast function to account for the correlation between successive increments.

Our analysis demonstrates that, despite the loss of Markovianity due to the velocity approximation in the partial observation case, the estimators converge to a Gaussian distribution (with a correction factor in the partial observation case). The proofs are based on Ito-like bounds and an adaptation of the Euler scheme for both the drift and diffusion components. Additionally, we provide insights into Hörmander’s condition, which helps establish hypoellipticity in our model within the framework of stochastic calculus of variations.
 

\end{abstract}

\tableofcontents

\section{Introduction}
In this paper, we consider the problem of parametric estimation for an interacting particle system of hypoelliptic diffusions. We focus on 
$N$ particles, each in dimension $2$, generalizing systems referred to as Langevin or hypoelliptic in various references in the literature.

Even without considering interacting particles and focusing solely on classical SDEs, such structures naturally arise in various applications across different domains, including finance, biology, and random mechanics. For instance, macromolecular systems can be modeled using these processes (see Grubmuller and Tavan (1994) \cite{GruTav} and Hummer (2005) \cite{Hum}). Another notable application is in audio signal analysis, as discussed by Giannopoulos and Godsill (2001) \cite{GiaGod}, with further methodological details provided by Godsill and Yang (2006) \cite{GodYan}.

Classical examples widely used in the literature include the stochastic growth model, the harmonic oscillator, and the oscillator with trigonometric potential. The latter describes the dynamics of a particle moving in a potential that is a superposition of trigonometric functions and is sometimes used in molecular dynamics in connection with the dynamics of dihedral angles (see Lasota and Mackey (1994) \cite{LasMac}). For a detailed discussion of such systems, we refer to Section 2 of \cite{PokStu}. 

One of the most classical applications is in kinetic systems, modeled by two-dimensional diffusion processes representing the position and velocity of an object. The natural extension of this is to consider kinetic systems of interacting particles, where $N$ such objects interact with each other, resulting in $N$ two-dimensional diffusions (see p. 44 in Section 3.1.3 of \cite{ChDi21} for details).

When Boltzmann published his seminal work \cite{24poc}, the study of large systems of interacting particles was primarily driven by the desire to model thermodynamic systems at a microscopic level. He argued that since a macroscopic volume of gas contains an immense number of elementary particles, tracking each individual particle is both impractical and unnecessary. Instead, a statistical approach should be adopted. Boltzmann, together with Gibbs, laid the foundation for a consistent kinetic theory of gases, which relies on a crucial assumption known as molecular chaos. This assumption posits that, despite the numerous interactions within the system, any two randomly selected particles should be statistically independent as the total number of particles approaches infinity. For an in-depth exploration of classical models in collisional kinetic theory, refer to Section 2.3.3 of \cite{Poc 1}.

Kac \cite{104poc} later introduced the concept that chaos, once established, should be propagated over time in evolving systems, a property known as the propagation of chaos. Shortly after Kac's work, McKean \cite{125poc} introduced a class of diffusion models that also satisfy this propagation of chaos property. The critical contribution of Kac and McKean was to demonstrate that the classical equations of kinetic theory have a natural stochastic interpretation. Their pioneering efforts have inspired the continued development of a rich and vibrant field in mathematical kinetic theory.

Over the past two decades and continuing today, systems of interacting particles have become ubiquitous across a wide range of applications. The tools and concepts originally developed in kinetic theory have transcended the confines of pure statistical physics. In fields such as mathematical biology and social sciences, self-organization models describe systems of indistinguishable particles—such as birds, insects, bacteria, and crowds—whose behavior is difficult to predict at the microscopic scale but can often be well-explained by models derived from the framework of mathematical kinetic theory (see, for example, \cite{2poc, 59poc, 132poc, 134poc}).

In another domain, the emerging theory of mean-field games explores the asymptotic properties of games involving large numbers of players \cite{30poc, 36poc}. More recently, particle systems have been applied to model complex phenomena in data science, with significant applications in neural network training \cite{47poc, 56poc, 142poc, 147poc}, optimization \cite{39poc, 90poc, 155poc}, and Markov Chain Monte Carlo methods \cite{65poc}. Notably, particle systems are increasingly employed to accelerate convergence in gradient-based methods or even to replace them entirely (see \cite{10f, 54f}).
For instance, in Section 2 of \cite{45f}, it is demonstrated that an algorithm based on a degenerate interacting particle system, known as particle gradient descent, can eliminate the need for interleaving Markov chains, while still providing provable guarantees, as discussed in \cite{13f}. Similarly, the recent work \cite{47f} (see Section 5) employs degenerate particle systems to develop an algorithm termed momentum particle descent. A closely related interacting particle system is introduced in \cite{2f}, where the authors propose the interacting particle Langevin algorithm, which is based on the discretization of interacting Langevin stochastic differential equations. Additionally, \cite{OlAk} investigates a kinetic interacting particle Langevin diffusion, corresponding to the underdamped setting of the overdamped Langevin diffusion proposed in \cite{2f}. {\rev General hypoelliptic models, such as the kinetic FitzHugh--Nagumo system, also arise naturally in neuroscience (see, for example,~\cite{MisQuiTou}).}

Given the vast range of applications associated with kinetic interacting particle systems, as extensively discussed above, we believe it is crucial to conduct a detailed statistical analysis of these equations, which is the primary focus of our work. In particular let us introduce, for $i \in \{1, ... , N \}$ and $t \in [0, T]$ the couple $Z_t^i := (X_t^i, Y_t^i) \in \R^2$ that satisfies
\begin{align}\label{eq: model intro}
\begin{cases}
d Y^i_t = b_1 (Z^i_t, \Pi^N_t ) dt,\\
d X^i_t = b_2 (Z^i_t, \Pi^N_t ) dt + a (Z^i_t, \Pi^N_t ) d B^i_t, 
\end{cases}
\end{align}
where $Z^1_0, \dots, Z^N_0$ are i.i.d.\ random vectors with a common law $\Pi_0$, independent of $(B^1_t)_{t \in [0,T]}, \dots, (B^N_t)_{t \in [0,T]}$, which in turn are independent standard Brownian motions, and 
$$
\Pi^N_t := \frac{1}{N} \sum_{i=1}^N \delta_{Z^i_t} \in \mathcal{P}_2 (\mathbb{R}^2)
$$
denotes the empirical measure of the system at time $t \in [0,T]$. The coefficients are such that $b_1, b_2, a: \R^2 \times \mathcal{P}_2(\R^2) \rightarrow \R$, where $\mathcal{P}_2(\R^2)$ denotes the set of probability measures on $\R^2$ with a finite second moment, endowed with the Wasserstein 2-metric
\begin{equation}
W_2(\mu, \nu) := \Big( \inf_{m \in \Gamma (\mu, \nu)} \int_{\R^2 \times \R^2} |x - y|^2 m(dx, dy) \Big)^{\frac 1 2},
\label{eq: wass}
\end{equation}
and $\Gamma(\mu, \nu)$ denotes the set of probability measures on $\R^2 \times \R^2$ with marginals $\mu$ and $\nu$. 

For \(N\) approaching infinity, the interacting particle system described above naturally converges to its mean field equation, represented by the McKean--Vlasov SDE:
\[
\begin{cases}
d \bar{Y}_t = b_1 (\bar{Z}_t, \bar{\Pi}_t ) \, dt,\\
d \bar{X}_t = b_2 (\bar{Z}_t, \bar{\Pi}_t ) \, dt + a (\bar{Z}_t, \bar{\Pi}_t ) \, dB_t, 
\end{cases} 
\]
for \(t \in [0,T]\), where \(\bar{Z}_t := (\bar{Y}_t, \bar{X}_t) \in \mathbb{R}^2\), and \(\bar{\Pi}_t\) denotes the law of \(\bar{Z}_t\). The initial condition \(\bar{Z}_0\) follows the distribution \(\bar{\Pi}_0 := \Pi_0\) and is independent of the standard Brownian motion \((B_t)_{t \in [0,T]}\).

This equation is nonlinear in the McKean sense (see e.g. \cite{125poc,79imp}), as the coefficients depend not only on the current state but also on the current distribution of the solution. It is well known that, under appropriate assumptions on the coefficients, a phenomenon known as \textit{propagation of chaos} can be observed, where the empirical law \(\Pi^N_t\) weakly converges to \(\bar{\Pi}_t\) as \(N \rightarrow \infty\) (see Section \ref{s: model} for further details and references). Such property will be crucial for our analysis.

We assume that the coefficients depend on some unknown parameters that we aim to estimate. Observe that, in the system \eqref{eq: model intro}, the noise directly affects the component \(X^i\), representing the 'speed' of the particle, and influences the position \(Y^i\) only through \(X^i\). It is important to note that in certain applications, it may not be possible to observe both the position and speed coordinates of each particle. This motivates us to study two different scenarios: in the first, we assume discrete observations of the paths of both \(X^i\) and \(Y^i\) for each particle; in the second, we consider only partial observations, where only the 'position' coordinates of the particles are observed.

It is worth emphasizing that, due to Hörmander's condition, we can introduce specific assumptions on the coefficients under which the model presented in \eqref{eq: model intro} becomes hypoelliptic in the context of stochastic calculus of variations. This allows us to transition to integrated diffusions (see Section \ref{s: hormander} for details), exemplified by the following system. For $i \in \{1, ... , N \}$ and $t \in [0, T]$ the couple $Z_t^i := (X_t^i, Y_t^i) \in \mathbb{R}^2$ satisfies
\begin{align}\label{eq:IPS intro}
\begin{cases}
d Y^i_t = X^i_t \, dt,\\
d X^i_t = b_{\mu_0} (Z^i_t, \Pi^N_t ) \, dt + a_{\sigma_0} (Z^i_t, \Pi^N_t ) \, d B^i_t,
\end{cases}
\end{align}
where $(\mu_0, \sigma_0)$ are the multidimensional parameters we aim to estimate. We refer to $Y^i_t$ and $X^i_t$ as the smooth and rough coordinates, respectively. {\spa While the theoretical reduction to \eqref{eq:IPS intro} is guaranteed by the implicit function theorem, the change of variables is generally not explicit. Consequently, estimating the parameters of the original system given observations of the integrated one (and vice versa) is a highly non-trivial task in practice. For this reason, we proceed by working directly under the assumption that the system takes the explicit integrated form \eqref{eq:IPS intro}. We stress that this is a standard structural assumption in the statistical inference of hypoelliptic diffusions, commonly required even in classical settings without interacting particles (cf.\ \cite[Assumption (C2)]{hypo}).} In the absence of interaction, diffusion processes of the form~\eqref{eq:IPS intro} are commonly known as \emph{integrated diffusions} (see, for instance, \cite{2hypo} or \cite{8hypo}). A notable instance of such models arises when the drift is given by $b(y,x) = - c(y,x)x - \nabla U(y)$ for some function $c$ and potential $U$, and the diffusion coefficient is constant. In this setting, the model~\eqref{eq:IPS intro} is known as a \emph{stochastic damping Hamiltonian system}, describing the motion of a particle under the combined influence of potential, dissipative, and random forces (see, e.g., \cite{Wu damp}).  
A classical example, in the absence of interaction, is the Kramers oscillator.  
Among hypoelliptic interacting diffusions, the most relevant example arising from applications and satisfying our hypothesis is the \emph{interacting underdamped Langevin dynamics}, describing a system of $N$ particles evolving in a confining potential, coupled through an interaction potential, and subject to both noise and dissipation (see Chapter~6 of \cite{47Igu}).

\subsection{Previous literature}
A large number of contribution is concerned with statistical inference for diffusions, we refer for example to the books \cite{Iac10, Kes12} and \cite{Kut04}. Moving in particular to parameter estimation for SDEs, a natural approach to estimating unknown parameters from the continuous observation of the diffusion would be to use maximum likelihood estimation. However, the likelihood function based on the discrete sample is not tractable, as it depends on the transition densities of the process, which are not explicitly known. {\spa To address this challenge, several methods have been developed for the high-frequency estimation of discretely observed classical SDEs. In this context, high-frequency refers to the asymptotic regime in which the discretization step $\Delta_n$ tends to zero as the number of observations $n$ increases.} A widely-used method involves considering a pseudo-likelihood function, often based on the high-frequency approximation of the process dynamics using the Euler scheme, as seen in \cite{Flo89, Kes97} and \cite{Yos92}.

Transitioning to statistical inference for interacting particle systems, it is worth noting that, aside from the early work by Kasonga in \cite{Kas}, the literature in this area is quite recent. The reason for this delay is that the initial interest in these systems stemmed from microscopic particle systems derived from statistical physics, which were not directly observable. However, the landscape has shifted due to the growing number of applications in this framework, as previously mentioned. These applications have generated data, sparking interest among statisticians to explore this area further.

Significant contributions to nonparametric or semiparametric statistical inference in this context can be found in works such as \cite{Pol, Vyt, Marc}, and \cite{Richard}. Meanwhile, the problem of parameter estimation from observations of both the interacting particle system and the associated McKean--Vlasov equation has been explored in various frameworks and asymptotic regimes; see, for instance, \cite{McKean, Bis, GenLar1, GenLar2, Giesecke, Liu, Sharrock, Wen} and the references therein.

 In all these works, the non-degeneracy of the diffusion coefficients plays a crucial role, since standard Euler-type contrast methods are not directly applicable in degenerate settings. To the best of our knowledge, at the time of our initial submission no results were available on statistical inference for degenerate interacting particle systems. During the revision process, however, we became aware of the recent and insightful works \cite{FonHof} and \cite{Igu new}, which investigate closely related frameworks.
In \cite{FonHof}, the authors analyze the empirical measure of a kinetic interacting particle system and establish a Bernstein-type inequality, which they exploit for the nonparametric estimation of the solution to the Vlasov–Fokker–Planck equation. Building on this, they further specialize to the FitzHugh–Nagumo model and derive optimal parameter estimators using moment-based methods. In contrast, \cite{Igu new} introduces a locally Gaussian approximation of the transition dynamics, carefully tailored to the degenerate structure of the diffusion matrix, thereby yielding a well-defined full likelihood. This enables parameter estimation for a broad class of hypoelliptic interacting particle systems.

{\rev Let us now review the existing literature in the absence of interaction. Ditlevsen and Sørensen~\cite{2hypo} studied discretely observed diffusion processes using prediction-based estimating functions, which were later applied to paleoclimate data in~\cite{DitDitAnd}. Gloter~\cite{8hypo} also addressed this problem, proposing a contrast function based on the Euler--Maruyama discretization. He suggested focusing on the rough coordinate of the SDE and approximating the unobserved integrated component via finite differences of the observed increments---an approach we will follow in this paper as well. This approximation introduces a bias, requiring a correction factor of $\frac{3}{2}$ in one of the terms of the contrast function under partial observations. This correction, however, increases the asymptotic variance of the diffusion parameter estimator. Samson and Thieullen~\cite{hypo} extended Gloter’s ideas to more general models, providing rigorous proofs. 

Pokern et al.~\cite{PokStu} proposed an empirical approximation of the likelihood based on the Itô--Taylor expansion, allowing the variance matrix to be invertible. However, their approach is limited to models where the drift is linear in the parameters and the volatility is constant. Ditlevsen and Samson~\cite{DitSam} focused on both filtering and inference for complete and partial observations, using an estimator based on the strong order 1.5 scheme from~\cite{KloPla}, which---similarly to~\cite{PokStu}---injects noise into the smooth coordinate. Notably,~\cite{DitSam} employed two separate contrast functions, leading to asymptotic marginal results rather than a joint central limit theorem.

Melnykova~\cite{Mel} later introduced a contrast estimator based on local linearization (LL), compared it with the least squares estimator, and established a joint central limit theorem. More recently, Gloter and Yoshida~\cite{GloYos, GloYos2} introduced adaptive and non-adaptive methods for hypoelliptic models, proving asymptotic normality for estimators based on complete observations. Their local Gaussian (LG) estimator generalizes the LL approach and builds upon~\cite{DitSam}, but with fewer restrictions on the model class. {\spa They demonstrated that employing a contrast function based only on the rough coordinate increases the asymptotic variance of the diffusion estimator, as the optimal variance is attained only when the full SDE is considered. A similar result was established in \cite{PilSamDit}.}

Further recent contributions include works by Iguchi et al.~\cite{Igu, Igu2, IguBes}, addressing limitations of earlier methods. In particular,~\cite{Igu2} introduced a closed-form contrast estimator based on Edgeworth expansions and Malliavin calculus;~\cite{Igu} focused on highly degenerate SDEs; and~\cite{IguBes} refined conditions for asymptotic normality, showing that the discretization step must tend to zero at a polynomial rate. Finally, the recent work~\cite{PilSamDit} proposes a novel estimator based on the Strang splitting scheme.}

\subsection{Our contribution}

{\rev Our objective is to construct estimators that are both consistent and asymptotically normal in both the partially and fully observed regimes, for the kinetic interacting particle system model defined in~\eqref{eq:IPS intro}.}

It is important to note that the presence of interaction introduces numerous novel aspects and challenges. Firstly, our asymptotic regime differs from the classical framework typically encountered in the SDE literature. {\spa In the classical regime, one usually considers a single process ($N=1$) over an increasing time horizon ($T \to \infty$), often requiring the existence of an invariant measure for consistency.} In our setting, the time horizon $T$ is fixed. Nonetheless, we are able to consistently estimate the drift and diffusion coefficients due to the increasing number of particles under consideration. Specifically, the number of particles $N$ plays a role analogous to that of $T$ in classical diffusion references, and similar conditions on the discretization step emerge, albeit involving $N$ instead of $T$ (see Remarks \ref{rk: conv rate} and \ref{rk: cond discretization} for details). Consequently, our approach does not rely on ergodicity to obtain the main results. Our asymptotics are derived as $N \to \infty$, crucially leveraging the propagation of chaos.

In this context, we propose two contrast functions tailored to the cases of complete and partial observations, respectively. Due to the non-invertibility of our diffusion matrix, the standard approach of constructing a quasi-likelihood based on the Euler--Maruyama discretization (as in \cite{McKean}) is not applicable. Nonetheless, our estimation procedure is still anchored in the Euler scheme, though applied differently. Specifically, we focus on estimating the drift and diffusion parameters of the second component of each particle, leading us to propose a contrast function solely based on the Euler approximation of this second component.

In the {\rev (most involving)} case of partial observations, where only \(Y^{i}\) is available, a correction term must be introduced in the contrast function. Our approach, detailed in Section \ref{s: main}, involves replacing the observation of \(X^i_{t_j}\) at each time \(t_j\) with the increment of the rate process \(\tilde{X}^i_{t_j} := \frac{1}{\Delta_n}(Y^i_{t_{j + 1}} - Y^i_{t_j})\), {\spa where $\Delta_n := t_{j+1} - t_j$ for $j = 0, \dots, n-1$ denotes the equidistant discretization step associated with a sample of size $n$.}

This substitution presents several challenges, notably the loss of Markovianity in the process $(\tilde{X}, Y)$. Additionally, the successive increments of the process \(\tilde{X}\) become correlated (see Remark \ref{rk: correlation tilde} for details). To address these issues, we modify the contrast function originally designed for the complete observation case, thereby defining the estimators for the drift and diffusion parameters.

To ensure the asymptotic normality of our estimators, we need to establish certain Ito-like bounds specific to our setting, which are outlined in Proposition \ref{prop: ito}. This poses significant challenges, not only because deriving these bounds requires differentiating with respect to the measure (necessitating an additional assumption—see Remark \ref{rk: assumption 7}), but also because applying the same method to the pair \((\tilde{X}, Y)\) introduces measurability issues related to the shifted definition of \(\tilde{X}\). Overcoming these difficulties required a deeper analysis of \(\tilde{X}\), which led us to consider higher-order Ito approximations, as described in Point (ii) of Proposition \ref{prop: ito}.

Through this approach—based on the analysis of functionals of \((X, Y)\) in the complete case and \((\tilde{X}, Y)\) in the partial case—we demonstrate that the proposed estimators for the drift and diffusion coefficients are consistent and asymptotically Gaussian in both scenarios.

The outline of the paper is as follows. In Section \ref{s: model}, we introduce the model and outline the assumptions necessary for our main results, along with an explanation of how the propagation of chaos applies in this context. Section \ref{s: main} presents our estimators and their properties, which constitute the core of our results. In Section \ref{s: hormander}, we discuss Hörmander's condition and explain the transition from model \eqref{eq: model intro} to \eqref{eq:IPS intro}, which is the model used for our estimation procedure. Section \ref{s: preparation} introduces several key results that are required for proving our main findings, while Section \ref{s: proof main} is dedicated to the proofs of these main results. Lastly, {\spa Appendix} \ref{s: proof technical} focuses on the detailed proofs of all the technical results introduced throughout the paper.

\textbf{Notation.} We write $E^\circ$ for the interior of a set $E \subset \R^d$. For $z := (z_1, \dots, z_d) \in \mathbb{R}^d$, we use the Euclidean norm $|z| := (z_1^2 + \dots + z_d^2)^{1/2}$. We denote by $\mathbb{R}^{2 \times 2}$ the set of 2-dimensional square real matrices.  For $A \in \mathbb{R}^{2 \times 2}$, we use the Frobenius norm $|A|:= (\sum_{i=1}^2 \sum_{j=1}^2 a_{ij}^2 )^{1/2}$. For $A,B \in \mathbb{R}^{2 \times 2}$, we have $|AB| \le |A||B|$. 
We denote by $\mathcal{P}_p (\mathbb{R}^2)$ the set of probability measures with bounded moment of order $p \ge 1$ on $\mathbb{R}^2$. Moreover, we use the Wasserstein $p$-distance between probability measures $\pi, \bar \pi \in \mathcal{P}_p (\mathbb{R}^2)$:
$$
W_p(\pi,\bar \pi) := \inf_{\gamma \in \Gamma(\pi,\bar \pi)} \Big( \int_{\mathbb{R}^2 \times \mathbb{R}^2} |z-\bar z|^p \gamma (d z, d \bar z) \Big)^{1/p}, 
$$
where $\Gamma (\pi, \bar \pi)$ denotes the set of probability measures on $\mathbb{R}^2 \times \mathbb{R}^2$ with the first, second marginals respectively equal to $\pi, \bar \pi$. \\
\\
\noindent Let us introduce $\mathcal{F}_t^N := \sigma \{(B_u^k)_{u \in [0,t]}, \, Z_0^{k}; \, k= 1, \dots , N \}$ and $\E_t [\cdot] := \E[\cdot | \mathcal{F}^N_t]$. For a set $(Y^{i,N}_{t,n})$ of random variables
and for a
$\delta \ge 0$, the notation 
$$
Y^{i,N}_{t,n} = R_t (\Delta_n^\delta)
$$ 
means that $Y^{i,N}_{t,n}$ is $\mathcal{F}_t^N$-measurable and the set $(Y^{i,N}_{t,n}/\Delta_n^\delta)$ is 
bounded in $L^q$ for all $q \ge 1$, uniformly in $t, i, n, N$, i.e.\ for all $q \ge 1$ there exists a $C_q>0$ such that for all $t, i, n, N$, 
$$
\E \big[ \big|Y^{i,N}_{t,n}/\Delta_n^\delta \big|^q \big]^{1/q}  \le C_q.
$$  
We point out that such notation is classical for the reminder term, in the context of parameter estimation for stochastic processes. \\
\\
We say that a function $f : \R^2 \times \mathcal{P}_2(\R^2) \to \R$ has a polynomial growth if there exist $C>0$, $k \in \mathbb{N}$ such that for all $(z,\pi) \in \R^2 \times \mathcal{P}_2(\R^2)$,
\begin{equation*}\label{eq: pol}
|f (z,\pi)| \le C (1+|z|^k +W_2^k(\pi,\delta_0)).
\end{equation*}
A function  $f : \R^2 \times \mathcal{P}_2(\R^2) \to \R$ satisfies a local Lipschitz condition and has a polynomial growth if there exist $C>0$, $k \in \mathbb{N}$ such that for all $(z_1,\pi_1), (z_2,\pi_2) \in \R^2 \times \mathcal{P}_2(\R^2)$,
\begin{align}
|f(z_1,\pi_1)-f(z_2,\pi_2)|&\le C (|z_1-z_2| + W_2(\pi_1,\pi_2))\nonumber\\ 
&\qquad\times (1 + |z_1|^k + |z_2|^k + W_2^k(\pi_1,\delta_0) + W_2^k(\pi_2,\delta_0)).\label{eq: cond pol growth}
\end{align}

\section{Model and assumptions}{\label{s: model}}
Let us begin by introducing the model we will be working with. As mentioned in the introduction, the model we aim to consider is the kinetic interacting particle system given in \eqref{eq: model intro}. However, in Section \ref{s: hormander}, we will elaborate on how, under certain hypotheses regarding the drift coefficients (see Assumption \ref{ass: hypo} in Section \ref{s: hormander}), this system is hypoelliptic in the sense of the stochastic calculus of variations, as discussed in detail in Section 2.3.2 of \cite{Nualart}. This property allows us to focus only on systems of interacting particles of the form described below, in \eqref{eq:IPS}.

Specifically, we assume that at time \( t \), the \( i \)-th particle is characterized by its position \( Y^i_t \) and its velocity \( X^i_t \). We further assume that the velocity is the time derivative of the position. More precisely, we assume the time evolution of the interacting particle system is given by
$$
Z^i_t :=  (Y^i_t,X^i_t) \in \mathbb{R}^2, 
\qquad i=1,\dots,N, \ t \in [0,T]
$$
and satisfies the system of SDEs:
\begin{align}\label{eq:IPS}
\begin{cases}
d Y^i_t = X^i_t dt,\\
d X^i_t = b (Z^i_t, \Pi^N_t ) dt + a (Z^i_t, \Pi^N_t ) d B^i_t, 
\end{cases}
\qquad i = 1,\dots, N, \ t \in [0,T].
\end{align}
We emphasize that, throughout our analysis, the time horizon $T$ will be fixed, and we will consider the asymptotic regime where only {\spa the number of particles $N$ and the sample size $n$ tend to infinity.}

A key point in our analysis is that we observe discrete observations of the interacting particles, but we can transition from the interacting particle system described above to a system of independent particles through what is known as the propagation of chaos. The system of independent particles consists of $N$ i.i.d. copies of the McKean--Vlasov equation associated with the model \eqref{eq:IPS}, as detailed below.

For $t \in [0,T]$, we consider a solution $\bar Z_t := (\bar Y_t, \bar X_t) \in \mathbb{R}^2$ of
\begin{align}\label{eq:MV}
\begin{cases}
d \bar Y_t = \bar X_t dt,\\
d \bar X_t = b (\bar Z_t, \bar \Pi_t ) dt + a (\bar Z_t, \bar \Pi_t ) d B_t, 
\end{cases}
\end{align}
where $\bar \Pi_t$ denotes the law of $\bar Z_t$, whereas $\bar Z_0$ has the law $\bar \Pi_0 := \Pi_0$ and is independent of the standard Brownian motion $(B_t)_{t \in [0,T]}$. 

Let us introduce the following assumptions. These are crucial to ensure that a solution to the SDE \eqref{eq:IPS} exists, possesses bounded moments of any order, and that the propagation of chaos applies.

\begin{assumption}{\label{as1}}
  $\Pi_0 \in \mathcal{P}_k (\mathbb{R}^2)$ for all $k \in \mathbb{N}$.   
\end{assumption}

\begin{assumption}{\label{as2}}
The functions
$$
b : \mathbb{R}^2 \times \mathcal{P}_2 (\mathbb{R}^2) \to \mathbb{R}, \qquad a : \mathbb{R}^2 \times \mathcal{P}_2(\mathbb{R}^2) \to \mathbb{R}
$$ 
satisfy the following Lischitz condition: there exists $C>0$ such that for all $(z,\pi), (\bar z,\bar \pi) \in \mathbb{R}^2 \times \mathcal{P}_2 (\mathbb{R}^2)$,
$$
\max (|b(z,\pi) - b(\bar z,\bar \pi)|, |a(z,\pi) - a(\bar z,\bar \pi)|) \le C (|z-\bar z| + W_2(\pi,\bar \pi)).
$$
\end{assumption}

\noindent Under the assumptions above, the propagation of chaos, as stated in \cite[Theorem 3.20]{ChDi21}, holds true. We restate it in Theorem \ref{th:PoC} below.

\begin{theorem}{\label{th:PoC}}
Assume \ref{as1}, \ref{as2}. Then, there exists $\varepsilon_N >0$ such that for all $i=1,\dots,N$, $t \in [0,T]$,
$$
\mathbb{E}[|Z^i_t-\bar Z^i_t|^2] \le \varepsilon_N
$$
and $\lim_{N \to \infty} \varepsilon_N = 0$. Here for $i=1,\dots,N$, $\bar Z^i := (\bar Z^i_t)_{t \in [0,T]}$ are independent copies of $\bar Z$ in \eqref{eq:MV}, starting from $\bar Z^i_0 := Z^i_0$ and driven by the same $B^i$ as $Z^i$.  
\end{theorem}
Note that \cite[Theorem 3.20]{ChDi21} applies because the drift and diffusion matrices of each particle, defined as:
$$
{\spa \tilde{B}} (z,\pi) := 
\begin{pmatrix}
x\\ 
b(z, \pi)    
\end{pmatrix} \in \mathbb{R}^2, \qquad {\spa \tilde{A}} (z, \pi) := 
\begin{pmatrix}
0 & 0\\
0 & a(z,\pi)
\end{pmatrix} \in \mathbb{R}^{2 \times 2}
$$
for all $z := (y,x) \in \mathbb{R}^2$, $\pi \in \mathcal{P}_2(\mathbb{R}^2)$, satisfy the Lipschitz condition.

{\spa We remark that the global Lipschitz assumption on the drift coefficient $b$ (Assumption A2) can be relaxed to accommodate locally Lipschitz drifts with polynomial growth. This relaxation is particularly relevant for interacting underdamped Langevin dynamics, where the potential function governing the drift is often non-linear, such as the standard double-well potential $V(x) = \frac{x^4}{4} - \frac{x^2}{2}$. In such a setting (cf.\ \cite[Assumption 2.1]{DosReis}), the boundedness of moments established in Lemma \ref{l: moments} can be recovered by applying \cite[Theorem 3.3]{AAP}, and the propagation of chaos follows from \cite[Proposition 3.1]{DosReis}. Consequently, the main results of this paper remain valid under these weakened conditions. 

Furthermore, throughout this work we focus on the one-dimensional kinetic setting where $Z_t^i = (X_t^i, Y_t^i) \in \mathbb{R}^2$, but the extension to the $d$-dimensional case ($Z_t^i \in \mathbb{R}^{2d}$) is theoretically consistent with our framework. Indeed, the propagation of chaos for such systems holds in any dimension, although the rate of convergence under the Wasserstein metric $W_p$ typically suffers from the curse of dimensionality, scaling as $N^{-1/d}$. However, since our derivation of the statistical rates depends on the convergence to the mean-field limit rather than the specific dimensional rate of that convergence, the results remain qualitatively unchanged. We have chosen to maintain the current one-dimensional notation to ensure the proofs remain accessible and to avoid the significant notational complexity associated with multidimensional matrix calculus.}

{\spa We recall that the coefficients depend on two parameters, $\mu_0$ and $\sigma_0$, which are the objects of our estimation.} Therefore, we will denote the drift coefficient as \(b_{\mu_0}\) and the diffusion coefficient as \(a_{\sigma_0}\). Furthermore, we will refer to the pair of parameters as \(\theta_0 := (\mu_0, \sigma_0)\), which belongs to the interior of a set \(\Theta := \Theta_1 \times \Theta_2\), where \(\Theta_i \subset \mathbb{R}^{p_i}\), \(i=1,2\), are compact and convex.


\section{Estimators and main results}{\label{s: main}}

Our aim is to estimate \(\theta_0\) from either complete or partial observations of the interacting particle system. In this section, we begin by detailing the two observation frameworks. In each framework, we will propose estimators that are specifically tailored to the discretization scheme under consideration. Finally, we will establish the consistency and asymptotic normality of the proposed estimators within their respective observation schemes.

We start by describing the complete observation regime, where we assume that for any \(i \in \{1, \dots, N\}\), both components are observed at discrete and equidistant times. Thus, we have access to
\[
(Y^i_{j \Delta_n}, X^i_{j \Delta_n}), \qquad i=1,\dots,N, \ j = 0, 1, \dots, n,
\]
where \(\Delta_n := T/n \to 0\) as \(n \to \infty\) and \(N \to \infty\), while \(T\) remains fixed.

The second case, referred to as the partial observation case, assumes that the processes \((X^i_t)_{t \in [0,T]}\) are hidden for all \(i \in \{1, \dots, N\}\). Therefore, we can only observe the positions
\[
Y^i_{j \Delta_n}, \qquad i=1,\dots,N, \ j = 0, 1, \dots, n.
\]
In this scenario, we need to pre-estimate the velocities \(X^i_{j \Delta_n}\) from the observed positions and use them in our estimation procedure. However, simply plugging these estimates into the contrast function used for complete observations would lead to a biased estimator. This necessitates the definition of a different contrast function for the partial observation case (see Remark \ref{rk: correlation tilde} for details).

This brings us to the definition of the contrast function within the two observation frameworks under consideration. {\spa Recall that in the idealized setting where a continuous-time path of the process is assumed to be available, the maximum likelihood estimator provides a highly efficient framework for parameter estimation.} However, in the realistic case where only a discrete-time sample of the process is observed, the transition density (and consequently the likelihood) is generally unavailable in closed form, making direct maximum likelihood estimation infeasible.

As pointed out in the introduction, a classical approach to overcoming this issue is to propose a quasi-likelihood function based on the Euler--Maruyama discretization of the process. Even in the case of complete observation, it is not feasible to directly apply the two-dimensional Euler contrast function for each particle \(i\) to estimate the parameters \((\mu_0, \sigma_0)\). Specifically, the two-dimensional Euler--Maruyama approximation for any \(i \in \{1, \dots, N\}\) is given by:
\[
\begin{pmatrix}
Y^i_{(j+1) \Delta_n}\\ 
X^i_{(j+1) \Delta_n}  
\end{pmatrix}
 = \begin{pmatrix}
Y^i_{j \Delta_n}\\ 
X^i_{j \Delta_n}  
\end{pmatrix}
+ \Delta_n {\spa \tilde{B}}(Z^i_{j \Delta_n}, \Pi^N_{j \Delta_n}) + \sqrt{\Delta_n} {\spa \tilde{A}}(Z^i_{j \Delta_n}, \Pi^N_{j \Delta_n}) \begin{pmatrix}
\xi^i_{j, 1}\\ 
\xi^i_{j, 2} 
\end{pmatrix}, 
\]
where $
\begin{pmatrix}
\xi^i_{j, 1}\\ 
\xi^i_{j, 2} 
\end{pmatrix}
$
is a vector of independent and identically distributed centered Gaussian variables.

However, the matrix \(A\) used above is not invertible, making this approach directly inapplicable. Nonetheless, our estimation procedure remains grounded in the Euler scheme. Rather than applying it directly, we focus on estimating the parameters in the drift and diffusion coefficients of the second component of each particle. Therefore, we propose a contrast function based solely on the Euler approximation of this second component. The remainder term arising from this approximation is analyzed in detail in Lemma 5.3 of \cite{McKean}. This leads us to the following contrast function:
\[ \label{eq: contrast complete}
\mathcal{L}^{N,C}_n (\theta) := \sum_{i=1}^N \sum_{j=0}^{n-1} \left( \frac{(X^i_{(j+1)\Delta_n} - X^i_{j\Delta_n} - \Delta_n b_\mu (Z^i_{j\Delta_n},  \Pi^N_{j\Delta_n}))^2}{\Delta_n a^2_\sigma (Z^i_{j\Delta_n}, \Pi^N_{j\Delta_n})} + \log \left( a^2_\sigma (Z^i_{j\Delta_n}, \Pi^N_{j\Delta_n}) \right) \right),
\]
where \(\theta = (\mu,\sigma)\) and the superscript \(C\) indicates {\spa ``complete observations''}. This function is an extension of the classical Euler contrast for interacting particle systems {\spa in one dimension}, as proposed in \cite{McKean}. Starting from such contrast function, we define a minimum contrast estimator $\hat{\theta}_n^{N,C}$ for complete observations as 
$$\hat{\theta}_n^{N,C} \in \mathop{\mathrm{arg\,min}}_{\theta \in \Theta} \mathcal{L}^{N,C}_n (\theta).$$

Let us now consider the case of partial observations, where the contrast function \eqref{eq: contrast complete} can no longer be applied due to the unavailability of the velocities \((X^i_{j \Delta_n})_{i=1, \dots, N; j =1, \dots, n}\). To address this challenge, various approaches have been explored in the literature. For instance, in \cite{8hypo}, Gloter proposes approximating \(X\) using the increments of \(Y\) in the context of integrated diffusion. This is the approach we have chosen to adopt, noting that in our framework, each particle also follows an integrated diffusion process.

To proceed, we introduce the increment (or rate) process:
$$
\tilde X^i_{j\Delta_n} := \frac{Y^i_{(j+1)\Delta_n}-Y^i_{j \Delta_n}}{\Delta_n}, \qquad \tilde Z^i_{j\Delta_n} := (Y^i_{j\Delta_n}, \tilde X^i_{j\Delta_n}) \qquad i = 1, \dots, N, \ j = 0, \dots, n-1,
$$
and
$$
\tilde \Pi^N_{j \Delta_n} := \frac{1}{N} \sum_{j=1}^N \delta_{\tilde Z^i_{j\Delta_n}}, \qquad j = 0, \dots, n-1.
$$

By definition, \(\tilde X^i_{j\Delta_n}\) depends solely on the observations of the positions. From the dynamics of the process \(Y^i\) described in \eqref{eq:IPS}, we have:
$$
\tilde X^i_{j\Delta_n} = \frac{1}{\Delta_n} \int_{j \Delta_n}^{(j + 1) \Delta_n} X_s^i \, ds,
$$
which justifies introducing \(\tilde X^i_{j\Delta_n}\) as a replacement for \(X^i_{j\Delta_n}\). For a small discretization step \(\Delta_n\), these two quantities are indeed close, as detailed in Proposition \ref{p: error X tilde}, which quantifies the error involved in this substitution. Despite the error being well-controlled, this substitution alone does not suffice to guarantee robust results when replacing \(X\) with \(\tilde{X}\) in the contrast function. As detailed in \cite{8hypo} for the classical diffusion case, this procedure would lead to a biased estimator. The reason is that successive terms of the rate process \(\tilde{X}\) are dependent (see Proposition \ref{prop: increments tilde} and Remark \ref{rk: correlation tilde} for details), necessitating a correction in the contrast function to account for this correlation. This brings us to the following contrast function:
\[ \label{eq: contrast partial}
\mathcal{L}^{N,P}_n (\theta) := \sum_{i=1}^N \sum_{j=1}^{n-2} \left( \frac{3}{2} \frac{(\tilde X^i_{(j+1)\Delta_n} - \tilde X^i_{j \Delta_n} - \Delta_n b_\mu (\tilde Z^i_{(j-1)\Delta_n}, \tilde \Pi^N_{(j-1)\Delta_n}) )^2}{\Delta_n a^2_\sigma (\tilde Z^i_{(j-1)\Delta_n}, \tilde \Pi^N_{(j-1)\Delta_n})} + \log (a^2_\sigma (\tilde Z^i_{(j-1)\Delta_n},\tilde \Pi^N_{(j-1)\Delta_n})) \right),
\]
where P stands for {\spa ``partial observations''}. Compared to the contrast function proposed for complete observations, the two key differences are the additional factor of \(\frac{3}{2}\), which corrects for the dependence structure highlighted above (and further detailed in Remark \ref{rk: correlation tilde}), and a shift in the index of the drift term.
It is important to note that \((\tilde X^i_{j\Delta_n})\) is not a Markov process, and the shift in the index of the drift coefficient has been introduced to avoid a correlation term between \((\tilde X^i_{(j+1)\Delta_n} - \tilde X^i_{j \Delta_n})\) and a function of \((\tilde Z^i_{j\Delta_n}, \tilde \Pi^N_{j\Delta_n})\). Although this correlation term is small (of the order \(\sqrt{\Delta_n}\)), it is not negligible in our analysis, leading us to the contrast function as defined above.

Finally, we define the minimum contrast estimator \(\hat{\theta}_n^{N,P}\) in the case of partial observation as follows:
$$
\hat{\theta}_n^{N,P} \in \mathop{\mathrm{arg\,min}}_{\theta \in \Theta} \mathcal{L}^{N,P}_n (\theta).
$$

With all this background we are ready to introduce some further assumptions, needed for our main results.

\begin{assumption} \label{as3}
(\textit{Regularity of the diffusion coefficient})
\textit{The diffusion coefficient is uniformly bounded away from $0$:
$$\inf_{(\sigma, z, \pi) \in \Theta_2 \times \R^2 \times {\cal P}_2(\R^2)} a^2_\sigma(z, \pi) >0. 
$$}
\end{assumption}
\begin{assumption} \label{as4}
(\textit{Regularity of the derivatives}) 
The first and second order derivatives in $\theta$ are locally Lipschitz in $(z,\pi)$ with polynomial weights, i.e. for all $\theta$ there exists $C >0$, $k,l =0,1,\dots$ such that for all $r_1 + r_2 =1, 2$, $h_1, h_2 = 1, ... , p_1$, $\tilde{h}_1, \tilde{h}_2 = 1, ... , p_2$,
$(z,\pi), (\bar{z},\bar{\pi}) \in \R^2 \times {\cal P}_2(\mathbb{R}^2)$,
\begin{align*}
&\big|\partial_{\mu_{h_1}}^{r_1} \partial_{\mu_{h_2}}^{r_2} b_\mu (z,\pi)- \partial_{\mu_{h_1}}^{r_1} \partial_{\mu_{h_2}}^{r_2} b_\mu (\bar{z},\bar{\pi}) \big| + \big|\partial_{\sigma_{\tilde{h}_1}}^{r_1} \partial_{\sigma_{\tilde{h}_2}}^{r_2} a_\sigma (z,\pi)-\partial_{\sigma_{\tilde{h}_1}}^{r_1} \partial_{\sigma_{\tilde{h}_2}}^{r_2} a_\sigma (\bar{z},\bar{\pi}) \big|
\\
&\qquad \le C (|z-\bar{z}| + W_2 (\pi,\bar{\pi})) \big( 1 + |z|^k + |\bar{z}|^k + W_2^l(\pi,\delta_0) + W_2^l(\bar{\pi},\delta_0) \big).
\end{align*}
\end{assumption}

Note that to establish the consistency of the proposed estimators in Theorem \ref{th: consistency} below, it would suffice to assume that the derivatives of the coefficients up to the third order exhibit polynomial growth, uniformly in the parameters. However, this assumption alone is insufficient to achieve the asymptotic normality of the estimators. For this, the stronger condition on the derivatives of the coefficients specified in Assumption \ref{as4} above is required.

Let us now introduce some notation that will be useful in the sequel. Recall that $\bar \Pi_t$ denotes the law of $\bar Z_t$ for $b=b_{\mu_0}$, $a = a_{\sigma_0}$. In the following we will write, for any function $f: \R^2 \times \mathcal{P}(\R^2) \rightarrow \R$ such that the mapping $(z,t) \mapsto f(z, \bar{\Pi}_t)$ is integrable with respect to $\bar{\Pi}_t (dz) \times dt$ over $\R^2 \times [0,T]$, 
$$\bar{\Pi}(f) := \int_0^{{T}} \Big( \int_{\R^2} f(z, { \bar{\Pi}_t}) \bar{\Pi}_t(d z) \Big) d t.$$
The integral $\bar \Pi(f)$ of a matrix-valued function $f$ is understood componentwise. We note that $(z,t) \mapsto f(z, \bar{\Pi}_t)$ is integrable if $f$ satisfies \eqref{eq: cond pol growth}. Indeed, \eqref{eq: cond pol growth} implies continuity of $(z,t) \mapsto f(z,\bar \Pi_t)$ and polynomial growth of $(z,\pi) \mapsto f(z,\pi)$, because $W_2(\bar \Pi_t, \bar \Pi_s) \le \E [|\bar Z_t - \bar Z_s|^2]^{\frac{1}{2}} \le C (t-s)^{\frac{1}{2}}$ and $\bar Z_t$ have bounded moments in $t$ by similar arguments as in the proof of our Lemma \ref{l: moments}(i), (ii) below.

We now state an assumption on the identifiability of the model. For this purpose we define the functions $I : \Theta \to \R$, $J : \Theta_2 \to \R$ as
\begin{align}\label{def:I}
I (\theta)&:= \bar{\Pi} \Big( \Big(\frac{b_{\mu}- b_{\mu_0}}{a_\sigma}\Big)^2 \Big),\\[1.5 ex] \label{def:J}
J (\sigma) &:= \bar{\Pi} \Big( \frac{a^2_{\sigma_0}}{a^2_{\sigma}} + \log a^2_{\sigma} \Big).
\end{align}

\begin{assumption} \label{as5}
 \textit{(Identifiability)} \textit{The functions $I, J$ defined above satisfy that for every $\varepsilon > 0$,
$$
\inf_{\theta \in \Theta: |\mu-\mu_{0}|\ge \varepsilon} I(\theta)>0\qquad \mbox{and } \inf_{\sigma \in \Theta_2: |\sigma-\sigma_{0}|\ge \varepsilon} (J(\sigma)-J(\sigma_0))>0.
$$}
\end{assumption}

Assumptions \ref{as1}--\ref{as5} are necessary to establish the consistency of our estimator as detailed in Theorem \ref{th: consistency}. While these assumptions are standard in the statistical analysis of stochastic processes, the identifiability condition in Assumption \ref{as5} warrants particular attention. The quantities $I(\theta)$ and $J(\sigma)$, central to this requirement, depend on the underlying law of the system, making them challenging to verify explicitly in practice.

{\spa This challenge is common in parameter estimation for interacting systems; for a thorough discussion, we refer to \cite[Section 2.4]{Hof2}, where the authors link global identifiability to the non-degeneracy of the Fisher information matrix \cite[Proposition 16]{Hof2}. To illustrate how these conditions can be checked, consider the example of an interacting underdamped Langevin dynamics with constant diffusion coefficient $a_\sigma (z, \pi) = \sqrt{2\gamma / \beta}$ and drift coefficient $b_\mu (z, \pi) = - \gamma x - V'(y) - \int_{\R \times \R} W'(y - v) \pi (d u, d v)$, where the potentials $V(y)$, $W(y)$ are known and $\mu = \gamma$, $\sigma = \sqrt{\gamma / \beta}$. Then
\begin{equation*}
I(\theta) = \frac{(\mu_0 - \mu)^2}{\sigma^2} \int_0^T \int_{\mathbb{R} \times \mathbb{R}} x^2 \bar{\Pi}_t(dx, dy) dt \qquad \mbox{and } J(\sigma) = \left( \frac{\sigma_0^2}{\sigma^2} + \log \sigma^2 \right) T.
\end{equation*}
Under Assumption \ref{as3}, $\inf_{\theta \in \Theta: |\mu-\mu_{0}|\ge \varepsilon} I(\theta) > 0$ simply requires that the system does not remain at the origin, i.e.\ the second moment of the particles remains positive over the time horizon. Moreover, $J(\sigma)$
satisfies $\inf_{|\sigma-\sigma_{0}|\ge \varepsilon} (J(\sigma)-J(\sigma_0))>0$ for all $T > 0$, since $\log x > 1 - \frac{1}{x}$ for all $x >0$, $x \neq 1$. Therefore, the identifiability conditions hold for this Langevin-type dynamics.  
Similarly, the conditions can be verified for the other examples in Section \ref{s: examples}, where the drift and diffusion coefficients depend linearly on the unknown parameters.}

\begin{theorem}{(Consistency)}
Assume \ref{as1}-\ref{as5}. Then the estimators are consistent in probability:
$$\hat{\theta}_n^{N,C} \xrightarrow{\mathbb{P}} \theta_0, \qquad \hat{\theta}_n^{N,P} \xrightarrow{\mathbb{P}} \theta_0,  \quad \mbox{as } n, N \rightarrow \infty.$$
\label{th: consistency}
\end{theorem}
In order to prove the asymptotic normality of our estimators we need some extra assumptions on the coefficients, as below.

\begin{assumption}{\label{as6}}
(\textit{Invertibility}) \textit{The $p_j \times p_j$ matrices
$\Sigma^{(j)}(\theta_0) = (\Sigma^{(j)}_{kl}(\theta_0))$, 
defined by
\begin{align*}
\Sigma^{(j)}_{kl} (\theta_0) {:}= \begin{cases}
\displaystyle \bar{\Pi} \Big( \frac{\partial_{\mu_{k}} b_{\mu_0} \partial_{\mu_{l}} b_{\mu_0}}{a^2_{\sigma_0}} \Big), \qquad &j=1, \, k, l = 1, \dots, p_1,
\\
\displaystyle \bar{\Pi} \Big( \frac{ \partial_{\sigma_{k}} a^2_{\sigma_0} \partial_{\sigma_{l}} a^2_{\sigma_0}}{a^4_{\sigma_0}} \Big) ,\qquad &j=2, \, k,l = 1,\dots, p_2,
\end{cases}
\end{align*}
satisfy
$\operatorname{det}(\Sigma^{{(j)}} (\theta_0)) \neq 0${, $j=1,2$.}
}
\end{assumption}

Assumption \ref{as6} imposes an invertibility condition that is essential for establishing asymptotic normality. It is important to note that, in this assumption, all coefficients are evaluated at the true parameter value. Additionally, we introduce a specific assumption regarding the form of the diffusion coefficient, which also applies solely at the true parameter value.

\begin{assumption} \label{as7}
 \textit{ (Integral condition on the diffusion coefficient) At $\sigma_0$, for all $(z,\pi)$ the diffusion coefficient takes the form 
 $$
 a_{\sigma_0}( z, \pi) {:} = \Tilde{a}
 \Big( z, \int_{\mathbb{R}^2} K(z, \bar{z}) \pi(d\bar{z})\Big)
 $$
for some functions $\tilde a, K \in C^2(\R^4; \R)$, whose derivatives of order $1$ and $2$ have polynomial growth. 
}
\end{assumption}

Assumption \ref{as7} is a technical requirement necessary to establish the final bound in Proposition \ref{prop: increments tilde}(i). This assumption is rooted in the application of an Itô formula specifically adapted for the framework under consideration, as detailed in Proposition \ref{prop: ito}. To apply this formula, it is crucial to differentiate with respect to the measure, in the second component of the diffusion coefficient. Assumption \ref{as7} is one way to address this challenge, as it allows differentiation with respect to the measure to be reduced to differentiating the function \( K \). However, this is not the only method available—there are various ways to differentiate a function of a measure, as discussed in the literature (see, e.g., \cite[Chapter 5]{36poc} or \cite[Chapter 2.2, Appendix A.1]{30poc}).

As explained in detail in Remark \ref{rk: assumption 7}, the situation simplifies when only the drift depends on an unknown parameter, allowing for a very general diffusion coefficient without the need for Assumption \ref{as7}.

\begin{theorem}{(Asymptotic normality)}
Assume \ref{as1}-\ref{as7}. Assume, moreover, that $N\Delta_n \rightarrow 0$ for $N,n \rightarrow \infty$. Then, in the complete observation case, 
$$(\sqrt{N}(\hat{\mu}_n^{N,C} - \mu_0), \sqrt{N
/ {\Delta_n}}(\hat{\sigma}_n^{N,C} - \sigma_0) ) \xrightarrow{d} \mathcal{N}(0, (\Sigma^{(1)}(\theta_0))^{-1} ) \otimes \mathcal{N}(0, 2 (\Sigma^{(2)}(\theta_0))^{-1}).$$
In the partial observation case, instead, 
$$(\sqrt{N}(\hat{\mu}_n^{N,P} - \mu_0),\sqrt{N
/ {\Delta_n}}(\hat{\sigma}_n^{N,P} - \sigma_0) ) \xrightarrow{d} \mathcal{N}(0, (\Sigma^{(1)}(\theta_0))^{-1}) \otimes \mathcal{N}(0, (9/4) (\Sigma^{(2)}(\theta_0))^{-1} ).$$
\label{th: as norm}
Here $\Sigma^{(j)}(\theta_0)$, $j = 1, 2$, are as defined in \ref{as6}.
\end{theorem}

\begin{remark}{\label{rk: proof}}
{\spa The proof of the main results is provided in Section \ref{s: proof main} and centers on analyzing functionals of \( Z^i_{t_j} = (Y_{t_j}^i, X_{t_j}^i) \) and \( \tilde{Z}_{t_j} = (Y_{t_j}^i, \tilde{X}_{t_j}^i) \) for \( i \in \{1, \ldots, N\} \) and \( j \in \{1, \ldots, n-1\} \). The asymptotic analysis of these functionals proves to be challenging due to the involvement of both the Markovian process \( Z^i_{t_j} \) and the non-Markovian process \( \tilde{Z}^i_{t_j} \).} 
\end{remark}


{\spa In our degenerate setting, the asymptotic variance involves the full law $\bar{\Pi}$ rather than merely the marginal of the second component. This mirrors the behavior of hypoelliptic diffusions in the ergodic $T \to \infty$ regime studied in \cite{hypo}, where the variance depends on the joint invariant density of $(X_t, Y_t)$. This contrasts with autonomous diffusions \cite{Kes97}, where the asymptotic distribution depends only on the stationary measure of the observed component. This parallelism underscores how our findings for $N \to \infty$ align with the hypoelliptic structure identified in \cite{hypo} (see Remark \ref{rk: conv rate} for further details).}


{\spa Our proposed drift estimator is efficient, as its asymptotic variance matches the Fisher information matrix derived from the Local Asymptotic Normality (LAN) property for $d$-dimensional McKean--Vlasov models under continuous observations \cite{Hof2}. Although the LAN property for joint estimation of drift and diffusion coefficients remains open, parallels with the classical $N=1$ and $T \to \infty$ setting (Remark \ref{rk: conv rate}) indicate that our diffusion estimator is inefficient. Specifically, constructing the contrast function solely from the rough coordinate inflates the asymptotic variance by a factor of $9/4$ relative to the optimal case. Thus, while consistent, the diffusion estimator lacks efficiency due to its partial use of the system's information.}

\begin{remark}{\label{rk: conv rate}}
It is well known that estimating the drift function from a single path of a stochastic differential equation over a fixed time horizon is not feasible. Therefore, when the goal is to estimate the drift, the time horizon \(T\) must tend to infinity to provide a sufficient amount of data. However, if observations from \(N\) independent paths are available, with \(N\) growing towards infinity, it becomes possible to estimate the drift over a fixed time horizon. 

As mentioned earlier, this is due to the parallelism between the roles played by \(T\) and \(N\) in the estimation procedure. Further evidence of this parallelism is found in the convergence rates presented in Theorem \ref{th: as norm}. Specifically, the convergence rate for drift estimation based on a single path over an infinite time horizon \(T = n \Delta_n\) is \(\sqrt{T}\), while the convergence rate for diffusion estimation is \(\sqrt{n} = \sqrt{T / \Delta_n}\), as detailed in Theorem 2 of \cite{hypo}. This clearly aligns with our results.
\end{remark}

\begin{remark}{\label{rk: cond discretization}}
{\spa The condition $N \Delta_n \to 0$ is essential for the asymptotic normality in Theorem \ref{th: as norm}, serving to approximate the contrast derivative with martingale increments and to manage inter-particle correlations. In the classical setting where $T = n\Delta_n \to \infty$, such constraints (e.g., $n\Delta_n^2 \to 0$) have been progressively relaxed using higher-order approximations \cite{Flo89, Yos92, Kes97, IguBes}, with similar developments for jump diffusions \cite{Sjs, Joint, GLM, Shi} and interacting particle systems \cite{McKean}. However, while higher-order schemes might mitigate discretization errors, they do not resolve the intrinsic particle correlations. In our partially observed framework, the dependency between $X$ and $\tilde{X}$ further complicates this relaxation, which we leave for future investigation.}
\end{remark}

{\spa Our primary goal is to initiate statistical inference for kinetic interacting particle systems, a field currently less mature than its non-interacting counterpart. In this work, we employ a locally degenerate approximation for the contrast function, based on the Euler--Maruyama discretization of the intractable process.

Alternative strategies in the non-interacting setting (see e.g.,~\cite{DitSam, GloYos, Igu, Mel, PilSamDit}) utilize locally non-degenerate approximations via transition density schemes. While such approaches accommodate broader hypoelliptic models, their asymptotic analysis for interacting particle systems is significantly more complex. To maintain a focused scope, we prioritize the locally degenerate framework. Notably, the very recent contribution~\cite{Igu new}, appearing during our revision, extends these ideas to the non-degenerate case. We view this as evidence of the timeliness of our work and its role in stimulating further research into broader classes of degenerate interacting particle systems.}

{\spa 
\section{Motivating model examples}\label{s: examples}

In this section, we provide several canonical examples from statistical physics, engineering, and econophysics to demonstrate the versatility of the model framework introduced in \eqref{eq:IPS intro}.

\subsection{Interacting underdamped Langevin dynamics}
The primary physical realization of our framework is the interacting underdamped Langevin dynamics. This model describes $N$ particles subject to a confining potential $V$, a reciprocal pair-interaction potential $W$, and a dissipative friction force. The drift and diffusion coefficients are given by:
\begin{equation}
    \begin{aligned}
        b(Z^i_t, \Pi^N_t) &= -\gamma X^i_t - V' (Y^i_t) - \frac{1}{N} \sum_{j=1}^N W' (Y^i_t - Y^j_t), \\
        a(Z^i_t, \Pi^N_t) &= \sqrt{2\gamma \beta^{-1}},
    \end{aligned}
\end{equation}
where $\gamma > 0$ is the friction coefficient and $\beta^{-1}$ represents the temperature of the environment. The drift term $b$ encapsulates the competition between conservative forces and non-conservative dissipation, serving as a microscopic foundation for Vlasov-type kinetic equations.

\subsection{The kinetic mean-field Ornstein--Uhlenbeck process}
A fundamental benchmark for analytical tractability and hypocoercivity theory is the ``harmonic oscillator'' of kinetic systems. By assuming quadratic potentials of the form $V(y) = \frac{\kappa}{2} y^2$ and $W(y) = \frac{\alpha}{2} y^2$, the drift becomes linear in the state variables:
\begin{equation} \label{eq:kinetic_OU}
    d X^i_t = \left( -\gamma X^i_t - \kappa Y^i_t - \alpha (Y^i_t - \bar{Y}^N_t) \right) dt + \sqrt{2\gamma \beta^{-1}} d B^i_t,
\end{equation}
where $\bar{Y}^N_t = \frac{1}{N} \sum_{j=1}^N Y^j_t$ represents the center of mass. Notice that $\bar{Y}^N_t$ is a direct function of the empirical measure $\Pi^N_t$, given by $\bar{Y}^N_t = \int_{\R \times \R} y \, \Pi^N_t (d x, d y)$. This model is unique as it preserves the Gaussianity of the distribution, providing a framework for studying convergence rates to equilibrium in degenerate diffusion settings.

\subsection{Inertial Kuramoto models in power systems}
Beyond statistical mechanics, the framework naturally extends to synchronization problems in power systems engineering. The synchronization of electrical generators can be mapped to our model by interpreting $Y^i$ as the phase angle and $X^i$ as the frequency deviation. Near the synchronized manifold, the sinusoidal coupling linearizes to:
\begin{equation}
    b(Z^i_t, \Pi^N_t) = -\gamma X^i_t - K (Y^i_t - \bar{Y}^N_t).
\end{equation}
This specific formulation is essential for assessing the resilience and stability of power grids against stochastic fluctuations in demand or renewable energy supply.

\subsection{Systemic risk in interbank networks}
Finally, a particular case of  \eqref{eq:kinetic_OU} is highly relevant in econophysics for modeling the stability of interbank reserves. Let $Y^i_t$ denote the log-monetary reserves of a given bank and $X^i_t$ its lending rate. Regulatory pressures and market averages act as harmonic restoring forces:
\begin{equation}
    b(Z^i_t, \Pi^N_t) = -\gamma X^i_t - \alpha (Y^i_t - \bar{Y}^N_t).
\end{equation}
Here, the ``harmonic'' attraction toward the market average $\bar{Y}^N_t$ represents a mean-reverting mechanism that prevents individual bank insolvency, while the noise $dB^i_t$ accounts for idiosyncratic market shocks.
}

\section{On the Hormander's condition for hypoellipticity}{\label{s: hormander}}
In this section, we aim to introduce specific assumptions on the coefficients under which the model presented in \eqref{eq: model intro} becomes hypoelliptic in the sense of the stochastic calculus of variations, as formally defined in Definition \ref{def: hypo} below. This will allow us to move to integrated diffusions, as in \eqref{eq:IPS}. To achieve this, some preliminary background is necessary. Our argumentation closely follows the approach proposed in Section 2.3.2 of \cite{Nualart}; readers interested in further insights on the topic may refer to that book.

{\spa First, let us introduce the empirical projection of a function defined on measures. For any real-valued function $\phi: \mathcal{P}_2(\mathbb{R}^l) \to \mathbb{R}$ (with $l \in \mathbb{N}$), we define its lifted version $\phi^N: \mathbb{R}^{l N} \to \mathbb{R}$ acting on the particle system by:
$$
\phi^N(z_1, \ldots, z_N) := \phi \left(\frac{1}{N} \sum_{k=1}^N \delta_{z_k}\right).
$$
This definition allows us to view the coefficients of the McKean--Vlasov equation as functions of the particle configuration $Z_t = (Z_t^1, \dots, Z_t^N)$. 

Specifically, for the drift coefficient $b_j(x, \mu)$ (where $j=1,2$), we define the corresponding function on the particle system, denoted by $b_j^{(i)}: \mathbb{R}^{2N} \to \mathbb{R}$, as:
$$
b_j^{(i)}(Z_t) := b_j\left(Z_t^i, \frac{1}{N} \sum_{k=1}^N \delta_{Z_t^k}\right).
$$
Here, the dependency on the measure argument is replaced by the empirical measure of the system $Z_t$. Similarly, the diffusion coefficient $a(x, \mu)$ induces a function $a^{(i)}: \mathbb{R}^{2N} \to \mathbb{R}$ defined by:
$$
a^{(i)}(Z_t) := a\left(Z_t^i, \frac{1}{N} \sum_{k=1}^N \delta_{Z_t^k}\right).
$$}
Thus, we can rewrite \eqref{eq: model intro} as the system given by the following $2N$ equations for $i \in \{1, \ldots, N \}$:

\begin{equation}{\label{eq: system 2-}}
\begin{cases}
d Y^i_t = b_1^{(i)} (Z_t) dt,\\
d X^i_t = b_2^{(i)} (Z_t) dt + a^{(i)} (Z_t) d B^i_t.
\end{cases}
\end{equation}
This system can also be represented in vector form as:
\begin{equation}\label{eq: system vect 2}
d Z_t = {\spa \hat{B}}(Z_t) dt + {\spa \hat{A}} (Z_t) d W_t,
\end{equation}
where $W_t$ is a $2N$-dimensional standard Brownian motion. Here, 
$${\spa \hat{B}}(Z_t) = (b_1^{(1)}(Z_t), b_2^{(1)}(Z_t), ... , b_1^{(N)}(Z_t), b_2^{(N)}(Z_t) ) \in \mathbb{R}^{2N},$$ 
and ${\spa \hat{A}} = ({\spa \hat{A}}_1 | {\spa \hat{A}}_2| ... | {\spa \hat{A}}_{2N})$ is a $2N \times 2N$ matrix where ${\spa \hat{A}_{2k + 1} (Z_t)} = \mathbf{0} \in \mathbb{R}^{2N}$ for $k \in \{0, ... , N-1 \}$, while for $k \in \{1, ... , N \}$, ${\spa \hat{A}_{2k} (Z_t)} = (0, ... , a^{(k)}(Z_t), ... , 0)^\top$, with the only nonzero value being $a^{(k)}(Z_t)$ at the $2k$-th coordinate. Recall that we are working under the hypothesis that $a(Z_t^k, \Pi_t^N)$ cannot be degenerate, which implies that $a^{(k)}(Z_t)$ is bounded away from $0$ for any $k \in \{1, ... , N \}$.

This notation is convenient because, to demonstrate the hypoellipticity of our system, we will impose non-degeneracy conditions on the coefficients in \eqref{eq: system vect 2}. These conditions will ensure the system is hypoelliptic according to Definition \ref{def: hypo} below. To introduce these constraints, consider the following vector fields on $\mathbb{R}^{2N}$ associated with the coefficients of Equation \eqref{eq: system vect 2}:
\begin{equation*}
\begin{cases}
A_j := \sum_{k = 1}^{2N} {\spa \hat{A}}_j^k (z) \partial_{z_k}, \quad j=1, ... , 2N, \\
B : = \sum_{k = 1}^{2N} {\spa \hat{B}}^k (z) \partial_{z_k}.
\end{cases}
\end{equation*}
The covariant derivative of $A_h$ in the direction of $A_j$ is defined as the vector field:
\begin{equation}\label{eq: cov der}
A_j^{\nabla} A_h := \sum_{l, k = 1}^{2N} {\spa \hat{A}}_j^l \partial_l {\spa \hat{A}}_h^k \partial_{z_k}.
\end{equation}
The Lie bracket between the vector fields $A_j$ and $A_h$ is defined by:
$$[A_j, A_h] := A_j^{\nabla} A_h - A_h^{\nabla} A_j.$$

We will now use these vector fields to express the stochastic differential equation \eqref{eq: system vect 2} in terms of the Stratonovich integral instead of the It\^{o} integral. For further details about the Stratonovich integral, we refer to \cite[Section 7]{ECinlar} and \cite{Thalmayer}. Given two real-valued continuous semimartingales $X$ and $Y$, the Stratonovich differential (see \cite[Definition 1.2.8]{Thalmayer}) is given by:
\begin{equation}
X \circ dY := X dY + \frac{1}{2} d[X,Y],
\end{equation}
where $dX$ is the standard It\^{o} differential and $[X,Y]$ is the quadratic covariation.

Under the Stratonovich formulation, the It\^{o} formula simplifies to the standard chain rule (see \cite[Proposition 1.2.10]{Thalmayer}). With this in mind, let us introduce the vector field $A_0$, which emerges from the It\^{o}-to-Stratonovich conversion of the model in \eqref{eq: system vect 2}:
\begin{equation}
A_0 := \sum_{k = 1}^{2N} \left( {\spa \hat{B}}^k(z) - \frac{1}{2} \sum_{j,l = 1}^{2N} {\spa \hat{A}}_l^j(z) \partial_j {\spa \hat{A}}_l^k(z) \right) \partial_{z_k} = B - \frac{1}{2} \sum_{l = 1}^{2N} A_l^\nabla A_l.
\end{equation}
It is straightforward to verify that rewriting \eqref{eq: system vect 2} as a Stratonovich equation yields:
\begin{equation}
dZ_t = A_0(Z_t) dt + \sum_{j = 1}^{2N} A_j(Z_t) \circ dW_t^j.
\end{equation}

At this point, we are equipped to define hypoellipticity in the context of the stochastic calculus of variations (see \cite[p. 129]{Nualart}).

\begin{definition}\label{def: hypo}
A differential operator $\mathcal{L}$ on an open set $G \subset \mathbb{R}^m$ with smooth coefficients is hypoelliptic if, for any distribution $\nu$ on $G$, $\nu$ is a smooth function on any open subset $G' \subset G$ where $\mathcal{L}\nu$ is smooth.
\end{definition}

Consider the second-order differential operator associated with our diffusion in \eqref{eq: system vect 2} for $m = 2N$:
\begin{equation}\label{eq: operator L}
 \mathcal{L} := \frac{1}{2} \sum_{j = 1}^{2N} (A_j)^2 + A_0. 
\end{equation}
H\"ormander's theorem (see \cite{138Nua}) states that if the Lie algebra generated by the vector fields $A_0, A_1, \ldots, A_{2N}$ has full rank at each point of $\mathbb{R}^{2N}$, then the operator $\mathcal{L}$ defined in \eqref{eq: operator L} is hypoelliptic.

Therefore, achieving the desired hypoellipticity of our model reduces to applying H\"ormander's theorem. This holds under a non-degeneracy hypothesis on the derivatives of the drift coefficients.

\begin{assumption}\label{ass: hypo}
For all $(y,x) \in \mathbb{R}^{N} \times \mathbb{R}^N$ and for all $k, i \in \{1, ... , N \}$, $\partial_{x_k} b_1^{(i)}(y,x) \neq 0$. 
\end{assumption}

{\spa
Although analyzing the hypoellipticity of individual particles as isolated two-dimensional SDEs is valid for the limiting McKean--Vlasov system (\cite[Section 2.1]{hypo}), this decoupled approach is unsuitable for interacting particles. Because the coefficients depend on the empirical measure, the noise across different particles is coupled, precluding a standard definition of the second-order differential operator for an isolated particle. Consequently, we investigate the hypoellipticity of the joint $2N$-dimensional SDE.

We also note that H\"ormander's condition in high dimensions has been established in related contexts. We refer to \cite[p. 366]{DitSam} and to the insightful discussions in \cite{Igu2, IguBes, Igu new}. However, our model differs structurally: in \cite{DitSam}, the rough component has dimension $p$ and the smooth component has dimension $1$, whereas here, both components of each particle have dimension $1$. This yields a $2N$-dimensional system where the rough and smooth components carry equal weight in the basis construction. Thus, while \cite{DitSam} builds a basis for $\mathbb{R}^{p+1}$, we require a basis for $\mathbb{R}^{2N}$. Although the underlying ideas are similar, this structural difference leads to technical distinctions (e.g., our Hypothesis A8 does not coincide with Condition 1 in \cite{DitSam}). We formalize our result in Proposition \ref{prop: hypo}, with the proof deferred to Section \ref{s: proof technical}.
}

\begin{proposition}\label{prop: hypo}
Assume that \cref{ass: hypo} holds true. Then, at every point $z \in \mathbb{R}^{2N}$, {\spa the vectors associated with the vector fields }
$$(A_2, [A_0, A_2], A_4, [A_0, A_4], \dots , A_{2N}, [A_0, A_{2N}]) =  (A_{2k}, [A_0, A_{2k}] )_{k \in \{1, \dots , N \}}$$
span the tangent space $T_z \mathbb{R}^{2N} \cong \mathbb{R}^{2N}$. Consequently, Hörmander's condition is satisfied and the system is hypoelliptic.
\end{proposition}

{\spa 
The algebraic condition established in Proposition \ref{prop: hypo} carries a clear probabilistic intuition that is well-known in the study of degenerate diffusions (see, e.g., \cite{Igu}). The fact that the vector fields $A_{2k}$ and their Lie brackets $[A_0, A_{2k}]$ span $\mathbb{R}^{2N}$ reflects how the stochastic driving forces propagate through the entire system. Specifically, the basis vectors $A_{2k}$ correspond to the direct Brownian noise $W^k_t$ acting on the velocity coordinates $X^k_t$. The Lie brackets $[A_0, A_{2k}]$, on the other hand, capture the integrated noise $\int_0^t W^k_s ds$. This integrated noise propagates into the position coordinates $Y^k_t$ through the kinematic relation $d Y^k_t = X^k_t dt$, as can be readily seen from an It\^{o}--Taylor expansion. Consequently, the noise effectively regularizes all $2N$ coordinates of the phase space, leading to the existence of a smooth Lebesgue density for the $2N$-dimensional system.
}

The proposition above is crucial for reducing the model \eqref{eq: system 2-} to an integrated diffusion. Specifically, one can apply the change of variable \(\hat{X}_t^{(i)} = b_1^{(i)}(Y_t, X_t)\). Consequently, for any \(i \in \{1, \ldots, N\}\), the first equation in \eqref{eq: system 2-} becomes:
\[dY_t^i = \hat{X}_t^i \, dt.\]
This suggests that \((Y_t^i, \hat{X}_t^i)\) should form an integrated diffusion. Then, thanks to \cref{ass: hypo}, one can apply the implicit function theorem, stating that \(X_t\) can be uniquely determined as a function of \((Y_t, \hat{X}_t)\), so that the vector \((Y_t^{(1)}, \hat{X}_t^{(1)}, \ldots, Y_t^{(N)}, \hat{X}_t^{(N)})\) satisfies
\begin{equation}\label{eq: model integrated}
\begin{cases}
dY_t^i = \hat{X}_t^i \, dt, \\
d\hat{X}_t^i = b(\hat{X}_t, Y_t) \, dt + a(\hat{X}_t, Y_t) \, dB_t^i,
\end{cases}
\end{equation}
for some coefficients \(b\) and \(a\), resulting from the implicit function theorem and Ito calculus. However, there is no general explicit expression available for \(X_t\) as a function of \((Y_t, \hat{X}_t)\). Therefore, one must assume that the original system in \eqref{eq: model intro} verifies a condition ensuring that, for \(i \in \{1, \ldots, N\}\), the process \((Y_t, \hat{X}_t)\), with \(\hat{X}_t^{(i)} = b_1^{(i)}(Y_t, X_t)\), satisfies the system in \eqref{eq: model integrated} with some explicit coefficients \(b\) and \(a\). Note that a similar condition is required even when considering a classical hypoelliptic diffusion, without the additional complexity of having \(N\) interacting particles (see Assumption (C2) in \cite{hypo}).

\section{Preparation for the proofs}{\label{s: preparation}}
Let us start by stating some moment inequalities that will be useful in the sequel. They are a direct consequence of Lemma 5.1 in \cite{McKean}. 
After that, we will analyze in detail what is the impact of replacing $X_{j \Delta_n}^i$ with $\tilde{X}_{j \Delta_n}^i$, in the partial observation case. In particular, we will deal with such error in Proposition \ref{p: error X tilde}, while Proposition \ref{prop: increments tilde} is devoted to the analysis of the increments of the process.

\begin{lemma}
Assume \ref{as1}-\ref{as2}. Let $p \ge 1$. Then for all $s, t \in [0,T]$ such that $t-s \in (0,1)$, $i = 1, ..., N$, $N \in \mathbb{N}$, the following 
hold true.
\begin{itemize}
    \item[(i)] 
    $\sup_{t \in [0,T]} \E[|Z_t^i|^p ] < C$,
    moreover, $\sup_{t \in [0,T]} \E[W_p^q(\Pi_t^{N}, \delta_0)] < C$ for $p \le q$.
    \item[(ii)]
    $\E[|Z_t^i - Z_s^i|^p] \le C (t - s)^{\frac{p}{2}}$.
    \item[(iii)] $\E[W_2^p(\Pi_t^{N}, \Pi_s^{N})] \le C (t - s)^{\frac{p}{2}}$.
\end{itemize}
\label{l: moments}
\end{lemma}




We now state a bound on the error committed moving from $X^i_{j \Delta_n}$ to $\tilde X^i_{j\Delta_n}$ and approximating $\tilde X^i_{j\Delta_n} - \tilde X^i_{(j-1)\Delta_n}$ by $\Delta_n b(\tilde Z^i_{j\Delta_n}, \tilde \Pi^N_{j\Delta_n})$.
Its proof is in Section \ref{s: proof technical}. For $i=1,\dots,N$, $j=0,\dots,n-1$,
let us introduce
\begin{align}{\label{eq: def xi tilde}}
\xi^i_j &:= \frac{1}{\Delta_n^{\frac 3 2}} \int_{j\Delta_n}^{(j+1)\Delta_n} ( (j+1)\Delta_n - s) d B^i_s,\\
{\label{eq: def xi}}
\tilde \xi^i_j &:= \frac{1}{\Delta_n^\frac{3}{2}} \int_{j \Delta_n}^{(j + 1)\Delta_n} (s - j \Delta_n) dB_s^i, \\
U^i_j &:= \tilde \xi^i_j + \xi^i_{j + 1}. \nonumber 
\end{align}

\begin{proposition}{\label{p: error X tilde}}
Assume \ref{as1}-\ref{as2}. If $k \ge 2$, then for all $i,j,N, n$, 
\begin{itemize}
    \item[(i)] $\tilde X^i_{j\Delta_n} - X^i_{j\Delta_n} = \Delta^{\frac{1}{2}}_n a_{\sigma_0} (Z^i_{j\Delta_n}, \Pi_{j \Delta_n}^N) \xi^i_j + \varepsilon^i_j$, where $\varepsilon^i_j$ is such that $\E [ |\varepsilon^i_j|^k] \le C \Delta_n^k$.
    \item[(ii)] $\E [ |\tilde X^i_{j\Delta_n} - X^i_{j\Delta_n}|^k] \le C \Delta_n^{\frac{k}{2}}$.
     \item[(iii)] 
     $\E [ W_2^k(\tilde{\Pi}_{j\Delta_n}^{N}, \Pi_{j\Delta_n}^{N})] \le C \Delta_n^{\frac{k}{2}}$.
\end{itemize}
\end{proposition}

\begin{remark}
Thanks to Proposition \ref{p: error X tilde}, we know that the error incurred by approximating the velocities with $\tilde{X}$ is of order $\Delta_n^\frac{1}{2}$. The final point, based on the same approximation, indicates that the error in approximating the empirical measure $\Pi_{j\Delta_n}^{N}$ (which depends on the unknown velocities) with $\tilde{\Pi}_{j\Delta_n}^{N}$ has the same magnitude.
\end{remark}


\begin{proposition}{\label{prop: increments tilde}}
 Assume \ref{as1}-\ref{as2}. If $k \ge 2$, then for all $i,j,N, n$,
 \begin{itemize}
     \item[(i)] $\tilde X^i_{(j+ 1)\Delta_n} - \tilde X^i_{j\Delta_n} - \Delta_n b_{\mu_0} (\tilde Z^i_{j\Delta_n}, \tilde{\Pi}_{j\Delta_n}^{N}) = \Delta_n^{\frac{1}{2}} a_{\sigma_0} (Z^i_{j\Delta_n},{\Pi}_{j\Delta_n}^{N}) U^i_j + \tilde{\varepsilon}_j^i,$ where $ \tilde{\varepsilon}_j^i$ is such that $\E [|\E_{j\Delta_n}[ \tilde{\varepsilon}_j^i]|^k] \le C \Delta_n^{\frac{3}{2} k}$ and $\E [ |\tilde{\varepsilon}_j^i|^k] \le C \Delta_n^k$. Assume moreover \ref{as7}. Then, $ \E [|\E_{j\Delta_n}[ \tilde{\varepsilon}_j^i U^i_j]|^k ] \le C \Delta_n^{\frac{3}{2}k}$.
     \item[(ii)] $\E [|\tilde X^i_{(j+ 1)\Delta_n} - \tilde X^i_{j\Delta_n}|^k] \le C \Delta_n^\frac{k}{2}$.
 \end{itemize}
\end{proposition}

\begin{remark}{\label{rk: correlation tilde}}
As a consequence of Proposition \ref{prop: increments tilde}, we observe that for any function of the two variables \(Y_t^i\) and \(X_t^i\), the term \(f(Z^i_{j \Delta_n}, {\Pi}_{j\Delta_n}^N)\) and \((\tilde{X}_{(j+1)\Delta_n}^i - \tilde{X}_{j\Delta_n}^i - \Delta_n b_{\mu_0}(\tilde{Z}_{j\Delta_n}^i, \tilde{\Pi}_{j\Delta_n}^N))\) have a correlation of order \(\Delta_n^\frac{1}{2}\). Although this correlation clearly tends to zero as \(n \to \infty\), it is not sufficiently small to not affect our estimation problem. This observation is what leads to the choice of the contrast function in \eqref{eq: contrast partial}. 

Additionally, it is important to note that this \(\Delta_n^\frac{1}{2}\) correlation arises from a stochastic integral, which can be approximated using the Euler--Maruyama scheme, as seen in the right hand side of (i). Crucially, this correlation term is conditionally centered, which plays a significant role in our analysis. It allows us to consider the following term in the development in (i), consisting in $\tilde{\varepsilon}_j^i$, which is now of order \(\Delta_n\).
\end{remark}

\begin{remark}{\label{rk: assumption 7}}
We would like to highlight that the assumption on the form of \(a\) stated in \ref{as7} is crucial for obtaining the final bound in (i). Specifically, it is necessary for proving certain Ito-like bounds tailored to our problem, which are detailed in Proposition \ref{prop: ito}. 

Without this additional hypothesis, it would be straightforward to verify from the proof that \([|\E_{j\Delta_n}[ \tilde{\varepsilon}_j^i U^i_j]|^k ] \le C \Delta_n^{k}\). However, this bound would not be sufficient to achieve the convergence in law required by Theorem \ref{th: convergence law QnN}, which is essential for establishing the asymptotic normality of \(\hat{\sigma}_n^{N,P}\).
It is then worth emphasizing that if one only seeks to estimate an unknown parameter in the drift component, no extra condition on the diffusion coefficient is necessary. In such cases, the diffusion's dependence on the measure can be as general as desired.
\end{remark}

\begin{proposition}{\label{prop: ito}}
Assume \ref{as1}-\ref{as2}. Let $f : \R^2 \times \mathcal{P}_2 (\R^2) \rightarrow \R$ satisfy \ref{as7} and assume its partial derivatives up to the second order have polynomial growth. Let $p \ge 1$. Then we have that for all $i,j,N, n$, 
\begin{itemize}
    \item[(i)] $\E[|\E_{(j - 1)\Delta_n}[f (Z^i_{j\Delta_n}, \Pi^N_{j\Delta_n})] - f (Z^i_{(j-1)\Delta_n}, \Pi^N_{(j-1)\Delta_n})|^p] \le C \Delta_n^p$,
    \item[(ii)] $\E[|\E_{(j - 1)\Delta_n}[f (\tilde{Z}^i_{(j-1)\Delta_n}, \tilde{\Pi}^N_{(j-1)\Delta_n})] - f (Z^i_{(j-1)\Delta_n}, \Pi^N_{(j-1)\Delta_n})|^p] \le C \Delta_n^p$.
\end{itemize}
\end{proposition}

The proof of our main results is based on the analysis of the asymptotic behaviour of some functionals we properly define below. Regarding the complete observation case let us introduce, for any function $f : \mathbb{R}^2 \times \mathcal{P}_2 (\mathbb{R}^2) \to \mathbb{R}$, the following: 
$$\nu_n^{N, C }(f) := \sum_{j = 0}^{n-1} \sum_{i = 1}^N f(Z_{j\Delta_n}^i, \Pi_{j\Delta_n}^N),$$
$$I_n^{N, C }(f) := \sum_{j = 0}^{n - 1} \sum_{i = 1}^N f(Z_{j\Delta_n}^i, \Pi_{j\Delta_n}^N)({X}_{(j+1)\Delta_n}^i - {X}_{j\Delta_n}^i - \Delta_n b_{\mu_0}(Z_{j \Delta_n}^i, \Pi_{j\Delta_n}^N)),$$
$$Q_n^{N, C }(f) := \sum_{j = 0}^{n - 1} \sum_{i = 1}^N f(Z_{j\Delta_n}^i, \Pi_{j\Delta_n}^N)({X}_{(j+1)\Delta_n}^i - {X}_{j\Delta_n}^i)^2.$$
Similarly, in order to study the case of partial observation, we introduce the following functionals:
$$\nu_n^{N, P}(f) := \sum_{j = 1}^{n -
2} \sum_{i = 1}^N f(\tilde{Z}_{(j-1)\Delta_n}^i, \tilde{\Pi}_{(j-1)\Delta_n}^N),$$
\begin{equation}{\label{eq: def InNP n+1}}
I_n^{N, P }(f) := \sum_{j = 1}^{n -2} \sum_{i = 1}^N f(\tilde{Z}_{(j-1)\Delta_n}^i, \tilde{\Pi}_{(j-1)\Delta_n}^N)(\tilde{X}_{(j+1)\Delta_n}^i - \tilde{X}_{j\Delta_n}^i - \Delta_n b_{\mu_0}(\tilde{Z}_{(j-1)\Delta_n}^i, \tilde{\Pi}_{(j-1)\Delta_n}^N)),
\end{equation}
$$Q_n^{N,P}(f) := \sum_{j = 1}^{n -2} \sum_{i = 1}^N f(\tilde Z_{(j-1)\Delta_n}^i, \tilde \Pi_{(j-1)\Delta_n}^N)(\tilde {X}_{(j+1)\Delta_n}^i - \tilde {X}_{j\Delta_n}^i)^2.$$

Note that in both \(I_n^{N, P }(f)\) and \(Q_n^{N,P}(f)\), we have introduced shifted processes within the function \(f\). This shift arises due to the correlation discussed in Remark \ref{rk: correlation tilde}, which led us to define the contrast function as in \eqref{eq: contrast partial}. Specifically, this definition implies that, when the square is expanded, certain functionals with the appropriate indices must be examined.

As highlighted in Remark \ref{rk: proof}, the analysis of these functionals is considerably more challenging compared to the case of classical (non-degenerate) interacting particles (as in \cite{McKean}) or the degenerate case with complete observation. This added complexity is due to the fact that we must now study functionals of \(\tilde Z_{(j-1)\Delta_n}^i\), which is non-Markovian, rather than the original Markovian process \(Z_{(j-1)\Delta_n}^i\).

Our main results are founded on the convergence of these functionals, as outlined in the following propositions. We begin by establishing the convergence in probability of the empirical mean \(\nu_n^N (f)\) in both the complete and partial observation settings. This result is heavily dependent on the propagation of chaos, which explains why the limit involves \(\bar{\Pi}_t\), the law of the process \(\bar{Z}_t\) evaluated at the true parameter values.

\begin{proposition}{\label{prop: converence nu}}
Assume \ref{as1}-\ref{as2}. Assume that for all $\theta \in \Theta$ the mapping $f_\theta : \R^2 \times \mathcal{P}_2 (\R^2) \rightarrow \R$ satisfies \eqref{eq: cond pol growth}. 
Moreover, assume that for all $(z,\pi)
$ the mapping $\theta \mapsto f_\theta (z,\pi)$ is continuously differentiable on $\Theta$
and that $(z,\pi) \mapsto \sup_{\theta \in \Theta} |\nabla_\theta f_\theta (z,\pi)|$ has a polynomial growth. 
Then we have that, uniformly in $\theta \in \Theta$, 
$$\frac{\Delta_n}{N} \nu_n^{N,C}(f_\theta) \xrightarrow{\mathbb{P}} \bar{\Pi}(f_\theta), \qquad \frac{\Delta_n}{N} \nu_n^{N,P}(f_\theta) \xrightarrow{\mathbb{P}} \bar{\Pi}(f_\theta).$$
\end{proposition}

Thanks to Proposition \ref{prop: converence nu}, we observe that although substituting the velocities with the approximating process \(\tilde{X}\) incurs an error (as detailed in Proposition \ref{p: error X tilde}), this error does not affect the limit of the empirical mean.

We now move to the analysis of $I_n^{N,C}(f)$ and $I_n^{N,P}(f)$, whose asymptotic behaviour is explained in the following theorem. 

\begin{theorem}{\label{th: convergence InN}}
Assume \ref{as1}-\ref{as2}. 
For all $(z,\pi) \in \R^2 \times \mathcal{P}_2(\R^2)$, let $\theta \mapsto f_\theta (z,\pi)$ be continuously differentiable on $\Theta$. Assume that for all $\theta \in \Theta$, 
$(z,\pi) \mapsto f_\theta (z,\pi)$ and $(z,\pi) \mapsto \sup_{\theta \in \Theta} |\nabla_\theta f_\theta (z,\pi)|$ have a polynomial growth. 
Then we have that, uniformly in $\theta \in \Theta$, 
$$\frac{1}{N} I_n^{N,C}(f_\theta) \xrightarrow{\mathbb{P}} 0, \qquad \frac{1}{N} I_n^{N,P}(f_\theta) \xrightarrow{\mathbb{P}} 0.$$  
\end{theorem}

Even in this case, the limit remains consistent in both scenarios of complete and partial observations, thanks to the shift introduced in the definition of \(I_n^{N,P}\). This shift allows us to circumvent the correlation term discussed in Remark \ref{rk: correlation tilde}. It is important to note that without this lag, the limit would not have been zero in the case of partial observations. This can be seen in \cite{8hypo}, where the author considers the hypoelliptic case without a shift for a single path of a stochastic differential equation. A careful reader will find that a similar approach applies to the interacting particle system, once the propagation of chaos is employed to remove the interactions among particles.

We now move to the analysis of the quadratic variations of $(Z_{t_j}^i)_{i,j}$ and $(\tilde{Z}_{t_j}^i)_{i,j}$. 

\begin{theorem}{\label{th: convergence QnN}}
Assume \ref{as1}-\ref{as2}. For all $\theta \in \Theta$, let $f_\theta: \R^2 \times \mathcal{P}_2 (\R^2) \rightarrow \R$ satisfy \eqref{eq: cond pol growth}. 
Moreover, assume that for all $(z,\pi)$ the mapping $\theta \mapsto f_\theta (z,\pi)$ is continuously differentiable on $\Theta$ and that $(z,\pi) \mapsto \sup_{\theta \in \Theta} |\nabla_\theta f_\theta (z,\pi)|$ has a polynomial growth. Then we have that, uniformly in $\theta \in \Theta$, 
$$\frac{1}{N} Q_n^{N,C}(f_\theta) \xrightarrow{\mathbb{P}} \bar{\Pi} ( f_\theta a^2_{\sigma_0} ), \qquad \frac{1}{N} Q_n^{N,P}(f_\theta) \xrightarrow{\mathbb{P}} \frac{2}{3} \bar{\Pi} ( f_\theta a^2_{\sigma_0} ).$$  
\end{theorem}

The results for complete observations extend to partial observations thanks to condition \eqref{eq: cond pol growth}. The approximation of $\tilde{X}^i_{(j+1)\Delta_n} - \tilde{X}^i_{j\Delta_n}$ in Proposition \ref{prop: increments tilde}(i) is also useful in the proofs of Theorems \ref{th: convergence InN} and \ref{th: convergence QnN}, for partial observations. Complete observations satisfy a similar approximation (see Lemma 5.3 and its proof in \cite{McKean}), but with $V^i_j := \Delta_n^{-1/2} (B^i_{(j+1)\Delta_n} - B^i_{j\Delta_n})$ replacing the error factor $U^i_j$. We note that $\E [(U^i_j)^2] = 2/3$, whereas $\E [(V^i_j)^2] = 1$. Therefore, we have an additional factor of $2/3$ in Theorem \ref{th: convergence QnN} for partial observations.

The convergence in probability established in the theorems above is essential for proving the consistency of our estimators. 
{\spa It is worth highlighting a direct and highly useful consequence of Theorem \ref{th: convergence QnN} under partial observations. When the target model features a constant volatility, that is, $a_\sigma(x) \equiv \sigma$, one can recover a quadratic variation type estimation of $\sigma$ by employing the process $\widetilde{X}$. Indeed, by setting $f_\theta \equiv 1$, Theorem \ref{th: convergence QnN} immediately yields the following convergence in probability:
$$ \frac{3}{2} \frac{1}{N} Q_n^{N, P}(1) \xrightarrow{\mathbb{P}} \sigma_0^2. $$}

We now turn to preparing the proof of asymptotic normality, which relies on the following convergence in distribution. We begin by analyzing \( I_n^{N,C}(f_\theta) \) and \( I_n^{N,P}(f_\theta) \).

\begin{theorem}{\label{th: convergence law InN}}
Assume \ref{as1}-\ref{as2}. Let $f : \R^2 \times \mathcal{P}_2 (\R^2) \rightarrow \R$ satisfy \eqref{eq: cond pol growth}.

Assume, moreover, that $N \Delta_n \rightarrow 0$. Then we have that
$$\frac{1}{\sqrt{N}} I_n^{N,C}(f) \xrightarrow{d} \mathcal{N}(0, \bar{\Pi} ( (f a_{\sigma_0})^2 )), \qquad \frac{1}{\sqrt{N}} I_n^{N,P}(f) \xrightarrow{d} \mathcal{N}(0, \bar{\Pi} ( (f a_{\sigma_0})^2)).$$  
\end{theorem}

In the theorems proving the convergence in distribution of the introduced functionals, we observe the condition \( N \Delta_n \to 0 \) appearing. It requires that the discretization step decreases rapidly enough. As discussed in Remark \ref{rk: cond discretization}, this condition is standard in the literature and is not specific to the degenerate case we are considering. However, it might be possible to relax this condition by using alternative (not explicit) contrast functions. For further details on this aspect, we refer again to Remark \ref{rk: cond discretization}.

We want to move to the convergence in distribution of $Q_n^{N,C}(f_\theta), Q_n^{N,P}(f_\theta)$. As we will see in the proof of Theorem \ref{th: convergence law QnN}, this heavily relies on the control gathered in Proposition \ref{prop: ito}, for which \ref{as7} is needed (see also Remark \ref{rk: assumption 7}).

\begin{theorem}{\label{th: convergence law QnN}}
Assume \ref{as1}-\ref{as2}. Let $f : \R^2 \times \mathcal{P}_2 (\R^2) \rightarrow \R$ satisfy \ref{as7} and \eqref{eq: cond pol growth}.
Assume, moreover, that 
$N \Delta_n \rightarrow 0$. Then we have that
\begin{align*}
&\frac{1}{\sqrt{N \Delta_n}} (Q_n^{N,C}(f) - \Delta_n \nu_n^{N,C}(f a_{\sigma_0}^2)) \xrightarrow{d} \mathcal{N}(0, 2 \bar{\Pi} ( f^2 a^4_{\sigma_0}) ),\\
& \frac{1}{\sqrt{N \Delta_n}} (Q_n^{N,P}(f) - \frac{2}{3} \Delta_n \nu_n^{N,P}(f a_{\sigma_0}^2)) \xrightarrow{d} \mathcal{N}(0, \bar{\Pi} ( f^2 a^4_{\sigma_0} )).
\end{align*}
\end{theorem}

From Theorem \ref{th: convergence QnN}, we see that in the case of partial observations, \( Q_n^{N, P}(f) \) underestimates \( \bar{\Pi} ( f_\theta a^2_{\sigma_0} ) \). This is why the central limit theorem includes the correction factor \(\frac{2}{3}\). A similar correction appears in Theorem 6 of \cite{hypo}.

\section{Proof {\spa of} main results}{\label{s: proof main}}
    Now that we have all the convergences of the functional we have introduced, as stated in previous section, we are ready to prove our main results. Let us start with the consistency gathered in Theorem \ref{th: consistency}.

\subsection{Proof of Theorem \ref{th: consistency}}
\begin{proof}
The proof follows closely the proof of the consistency in \cite{McKean}. 
We have to show that, uniformly in $(\mu, \sigma) = \theta \in \Theta$, 
\begin{align}
\frac{\Delta_n}{N} L_n^{N,C}(\mu,\sigma) & \xrightarrow{\mathbb{P}} \bar{\Pi} \Big( \frac{a^2_{\sigma_0}}{a^2_\sigma} + \log a^2_\sigma \Big), \label{eq: consistency sigma complete} \\
\frac{{\spa \Delta_n}}{N} L_n^{N,P}(\mu,\sigma) & \xrightarrow{\mathbb{P}} \bar{\Pi} \Big( \frac{a^2_{\sigma_0}}{a^2_\sigma} + \log a^2_\sigma \Big) \label{eq: consistency sigma partial}.
\end{align}
This guarantees the convergence of $\hat{\sigma}^2_n$ to $\sigma_0^2$ in both cases of complete and partial observations. Regarding the consistency of $\hat{\mu}_n$, this is ensured by the uniform in $(\mu,\sigma)$ convergence:
\begin{align}
\frac{1}{N} (L_n^{N,C}(\mu, \sigma) - L_n^{N,C}(\mu_0, \sigma) ) & \xrightarrow{\mathbb{P}} \bar{\Pi} \Big( \frac{(b_\mu - b_{\mu_0})^2}{a^2_\sigma} \Big) \label{eq: consistency mu complete}, \\
\frac{1}{N} (L_n^{N,P}(\mu, \sigma) - L_n^{N,P}(\mu_0, \sigma) ) & \xrightarrow{\mathbb{P}} \frac{3}{2} \bar{\Pi} \Big( \frac{(b_\mu - b_{\mu_0})^2}{a^2_\sigma} \Big) \label{eq: consistency mu partial}.
\end{align}

We start by showing \eqref{eq: consistency sigma complete}-\eqref{eq: consistency sigma partial} and so the consistency of $\hat{\sigma}_n^2$. Rewriting $L^{N,C}_n (\mu,\sigma)$ we see that
\begin{align}
\frac{\Delta_n}{N} L_n^{N,C}(\mu,\sigma) 
&= \frac{1}{N} Q_n^{N,C} \Big( \frac{1}{a^2_\sigma} \Big) - 2 \frac{\Delta_n }{N} I_n^{N,C} \Big( \frac{b_\mu}{a^2_\sigma} \Big) + \frac{\Delta_n^2}{N} \nu^{N,C}_n \Big( \frac{b_\mu}{a^2_\sigma} (b_\mu - 2 b_{\mu_0}) \Big) \nonumber\\ 
&\qquad+  \frac{\Delta_n}{N} \nu^{N,C}_n ( \log a^2_\sigma )
 \xrightarrow{\mathbb{P}} 
{\bar \Pi} \Big( \frac{a^2_{\sigma_0}}{a^2_\sigma} + \log a^2_\sigma \Big)
{\label{eq: start cons 1}}
\end{align}
uniformly in $\theta$ according to Proposition \ref{prop: converence nu} and Theorems \ref{th: convergence InN} and \ref{th: convergence QnN}.
Similarly, in the partial observation case we have that
\begin{align}{\label{eq: start partial 2}}
\frac{\Delta_n}{N} L_n^{N,P}(\mu,\sigma) &= \frac{3}{2} \frac{1}{N} Q_n^{N,P} \Big( \frac{1}{a^2_\sigma} \Big) - 3 \frac{\Delta_n}{N} I_n^{N,P} \Big( \frac{b_\mu}{a^2_\sigma} \Big) + \frac{3}{2} \frac{\Delta_n^2}{N} \nu_n^{N,P} \Big( \frac{b_\mu}{a^2_\sigma}(b_\mu - 2 b_{\mu_0}) \Big)\\ 
&\qquad+ \frac{\Delta_n}{N} \nu_n^{N,P}(\log a^2_\sigma ) \xrightarrow{\P} \bar{\Pi} \Big( \frac{a^2_{\sigma_0}}{a^2_\sigma} + \log a^2_\sigma \Big)\nonumber
\end{align}
uniformly in $(\mu,\sigma)$ according to
the same Proposition \ref{prop: converence nu} and Theorems \ref{th: convergence InN} and \ref{th: convergence QnN}.

Let us move to the proof of the consistency of $\hat{\mu}_n^N$. Observe that, according to the decomposition in \eqref{eq: start cons 1}, it is 
\begin{align*}
\frac{1}{N}(L_n^{N,C}(\mu, \sigma) - L_n^{N,C}(\mu_0, \sigma)) &= \frac{2}{N} I_n^{N,C} \Big( \frac{b_{\mu_0} - b_\mu}{a^2_\sigma} \Big) + \frac{\Delta_n}{N} \nu_n^{N,C} \Big(\frac{(b_\mu - b_{\mu_0} )^2}{a^2_\sigma} \Big)\\
&\xrightarrow{\P} \bar{\Pi} \Big(\frac{(b_\mu - b_{\mu_0} )^2}{a^2_\sigma} \Big)
\end{align*}
uniformly in $(\mu,\sigma)$ thanks to Proposition \ref{prop: converence nu} and Theorem \ref{th: convergence InN}.
Regarding the partial observation case, from \eqref{eq: start partial 2} we obtain 
\begin{align*}
\frac{1}{N}(L_n^{N,P}(\mu, \sigma) - L_n^{N,P}(\mu_0, \sigma)) &= \frac{3}{N} I_n^{N,P}\Big( \frac{b_{\mu_0} - b_\mu}{a^2_\sigma} \Big) + \frac{3}{2}\frac{\Delta}{N} \nu_n^{N,P}\Big( \frac{(b_\mu - b_{\mu_0} )^2}{a^2_\sigma} \Big)\\
&\xrightarrow{\P} \frac{3}{2}\bar{\Pi} \Big( \frac{(b_\mu - b_{\mu_0} )^2}{a^2_\sigma} \Big)
\end{align*}
uniformly in $(\mu,\sigma)$ again according to
Proposition \ref{prop: converence nu} and Theorem \ref{th: convergence InN}. 
It leads us to the consistency of $\hat{\mu}_n^N$ towards $\mu_0$, acting as in the proof of Theorem 3.1 in \cite{McKean}.

\end{proof}

\subsection{Proof of Theorem \ref{th: as norm}}
\begin{proof}
{\spa In what follows, we present the proofs of our main results for the case of scalar parameters in order to simplify the notation. For the same reason, we will omit in $\Delta_n$ the dependence in $n$.} The scheme of the proof is the same in both the complete and the partial observation case. For $A \in \{ P, C \}$, we use Taylor's formula and, as $\nabla_\theta L_n^{N,A}(\hat{\theta}_n^{N,A}) = 0$, we obtain 
$$(\hat{\theta}_n^{N,A} - \theta_0) \int_0^1 \nabla_\theta^2 L_n^{N,A}(\theta_0 + u ( \hat{\theta}_n^{N,A}- \theta_0)) du = - \nabla_\theta L_n^{N,A}(\theta_0). $$
Then, let us introduce the matrix 
\begin{equation*}
    M_n^N := \begin{pmatrix}
        \frac{1}{\sqrt{N}} & 0 \\
        0 & \sqrt{\frac{\Delta}{N}}
    \end{pmatrix}.
\end{equation*}
Multiplying the equation above by $M_n^N$ we obtain the equation 
\begin{equation}{\label{eq: start as norm}}
(\hat{\theta}_n^{N,A} - \theta_0) (M_n^N)^{-1} \int_0^1 \Sigma_n^{N,A}(\theta_0 + u ( \hat{\theta}_n^{N,A}- \theta_0)) du = - \nabla_\theta L_n^{N,A}(\theta_0) M_n^{N,A},
\end{equation}
where 
\begin{align*}
 \Sigma_n^{N,A}(\theta) & = M_n^{N,A}  \nabla^2_\theta L_n^{N,A}(\theta) M_n^{N,A} \\
 & = \begin{pmatrix}
     \frac{1}{N} \partial_\mu^2 L_n^{N,A}(\theta) & \frac{\sqrt{\Delta}}{N} \partial_\mu \partial_\sigma L_n^{N,A}(\theta) \\
     \frac{\sqrt{\Delta}}{N} \partial_\mu \partial_\sigma L_n^{N,A}(\theta) & \frac{{\Delta}}{N} \partial^2_\sigma L_n^{N,A}(\theta)
 \end{pmatrix}.
\end{align*}
The asymptotic behavior of the estimator is given by the asymptotic behavior of the first and second derivatives of the contrast function. Let us now detail the cases of complete and partial observation separately. We start by considering the complete observation framework. Note that 
\begin{equation*}
\nabla_\theta L_n^{N,C}(\theta_0) M_n^{N,C} = \begin{pmatrix}
\frac{1}{\sqrt{N}} \partial_\mu L_n^{N,C}(\theta_0) & \sqrt{\frac{\Delta}{N} }\partial_\sigma L_n^{N,C}(\theta_0)   
\end{pmatrix}.
\end{equation*}
According to \eqref{eq: start cons 1} it is easy to check that 
$$- \frac{1}{\sqrt{N}} \partial_\mu L_n^{N,C}(\theta_0) = \frac{2}{\sqrt{N}} I_n^{N,C}\Big( \frac{\partial_\mu b_{\mu_0}}{a^2_{\sigma_0}} \Big)$$
and
\begin{align*}
- \sqrt{\frac{\Delta}{N} }\partial_\sigma L_n^{N,C}(\theta_0) 
& = \frac{1}{\sqrt{N \Delta}} Q_n^{N,C} \Big(\frac{\partial_\sigma a^2_{\sigma_0}}{a^4_{\sigma_0}}\Big) - \sqrt{\frac{\Delta}{N}} \nu_n^{N,C} \Big(\frac{\partial_\sigma a^2_{\sigma_0}}{a^2_{\sigma_0}}\Big) + o_{\mathbb{P}}(1).
\end{align*}
{\spa In the equality above,  $o_{\mathbb{P}}(1)$ denotes the terms 
\begin{equation*}
-2 \sqrt{\frac{\Delta}{N}} I_n^{N,C} \Big( \frac{ \partial_\sigma a^2_{\sigma_0} b_{\mu_0}}{a^4_{\sigma_0}} \Big) - \Delta \sqrt{\frac{\Delta}{N}} \nu_n^{N,C} \Big(\frac{ \partial_\sigma a^2_{\sigma_0} b_{\mu_0}^2 }{a^4_{\sigma_0}}\Big),
\end{equation*}
both of which vanish in probability as $n, N \to \infty$. This negligibility is a direct consequence of the assumption $N \Delta \to 0$, as well as the convergence results established in Theorem \ref{th: convergence InN} and Proposition \ref{prop: converence nu}, respectively.}
Then, Theorem \ref{th: convergence law InN} ensures that 
$$- \frac{1}{\sqrt{N}} \partial_\mu L_n^{N,C}(\theta_0) \xrightarrow{d} \mathcal{N} \Big(0, \bar{\Pi} \Big( \Big(2 \frac{\partial_\mu b_{\mu_0}}{(a^2_{\sigma_0})}a_{\sigma_0} \Big)^2 \Big) \Big) = N \Big(0,4 \bar{\Pi}\Big( \Big(\frac{\partial_\mu b_{\mu_0}}{a_{\sigma_0}}\Big)^2 \Big) \Big). $$
Moreover, Theorem \ref{th: convergence law QnN} implies 
$$- \sqrt{\frac{\Delta}{N} }\partial_\sigma L_n^{N,C}(\theta_0) \xrightarrow{d} \mathcal{N} \Big(0, 2 \bar{\Pi} \Big( \Big( \frac{\partial_\sigma a^2_{\sigma_0}}{a^2_{\sigma_0}} \Big)^2 \Big) \Big). $$
We deduce,
$$\nabla_\theta L_n^{N,C}(\theta_0) M_n^{N,C} \xrightarrow{d} \mathcal{N}(0, \Sigma^C(\theta_0)),$$
with 
\begin{equation*}
 \Sigma^C(\theta_0) = \begin{pmatrix}
     4 \bar{\Pi}\big( \big(\frac{\partial_\mu b_{\mu_0}}{a_{\sigma_0}}\big)^2\big) & 0 \\
     0 & 2 \bar{\Pi}\big(\big(\frac{\partial_\sigma a^2_{\sigma_0}}{a^2_{\sigma_0}}\big)^2\big)
 \end{pmatrix}.
\end{equation*}
Observe that the joint asymptotic normality of $\partial_\mu L^{N,C}_n (\theta_0)$ and $\partial_\sigma L^{N,C}_n (\theta_0)$ comes from their asymptotic independence. Let us move to the analysis of the second derivatives. We are going to prove that $\Sigma_n^{N,C}(\theta)$ converges in probability, uniformly in $\theta$, towards $\Sigma(\theta)$. 
Indeed, 
\begin{align*}
\frac{1}{N} \partial^2_\mu L_n^{N,C}(\theta) &= - \frac{2}{N} I_n^{N,C} \Big(\frac{\partial_\mu^2 b_\mu}{a_\sigma^2}\Big) + \frac{\Delta}{N} \nu_n^{N,C} \Big(\partial_\mu \big(2 \frac{\partial_\mu b_\mu}{a_\sigma^2}(b_\mu - b_{\mu_0})\big)\Big) \\
& = o_{\mathbb{P}}(1) + \frac{\Delta}{N} \nu_n^{N,C}\Big(\frac{2}{a_\sigma^2}\big(\partial_\mu^2 b_\mu(b_\mu - b_{\mu_0}) + (\partial_\mu b_\mu)^2\big)\Big) \\
& \xrightarrow{\mathbb{P}} \bar{\Pi}\Big(\frac{2}{a_\sigma^2}\big(\partial_\mu^2 b_\mu(b_\mu - b_{\mu_0}) + (\partial_\mu b_\mu)^2\big)\Big) =: 2 \Sigma^{(1)}(\theta),
\end{align*}
{\spa where the $o_{\mathbb{P}}(1)$ term arises because $- \frac{2}{N} I_n^{N,C} (\frac{\partial_\mu^2 b_\mu}{a_\sigma^2})$ converges to zero in probability, an asymptotic property guaranteed by Theorem \ref{th: convergence InN}.} The convergence holds true uniformly in $\theta$ thanks to Proposition \ref{prop: converence nu}. Observe in particular that $\Sigma^{(1)}(\theta_0) = \bar{\Pi}(\frac{(\partial_\mu b_{\mu_0})^2}{a^2_{\sigma_0}} )$. Furthermore, 
\begin{align*}
\frac{\sqrt{\Delta}}{N} \partial_\mu \partial_\sigma L_n^{N,C}(\theta) & = \frac{\sqrt{\Delta}}{N} \partial_\sigma \partial_\mu L_n^{N,C}(\theta) \\
& = 2 \frac{\sqrt{\Delta}}{N} I_n^{N,C} \Big( \frac{\partial_\mu b_\mu}{a_\sigma^4} \partial_\sigma a_\sigma^2 \Big) - \frac{\Delta^\frac{3}{2}}{N} \nu_n^{N,C} \Big( 2 \frac{\partial_\sigma a_\sigma^2 \partial_\mu b_\mu}{a_\sigma^4}(b_\mu - b_{\mu_0}) \Big) \xrightarrow{\mathbb{P}} 0
\end{align*}
uniformly in $\theta$ as a consequence of Proposition \ref{prop: converence nu} and Theorem \ref{th: convergence InN}. \\
Regarding the second derivatives in $\sigma$, 
\begin{align*}
\frac{\Delta}{N} \partial^2_\sigma L_n^{N,C}(\theta) = \frac{1}{N} Q_n^{N,C} \Big( 2 \frac{(\partial_\sigma a_\sigma^2)^2}{a_\sigma^6} - \frac{\partial^2_\sigma a_\sigma^2}{a^4_\sigma} \Big) + o_{\mathbb{P}}(1) + \frac{\Delta}{N} \nu_n^{N,C} \Big( \frac{\partial^2_\sigma a_\sigma^2}{a_\sigma^2}- \frac{(\partial_\sigma a_\sigma^2)^2}{a_\sigma^4} \Big),
\end{align*}
{\spa where we have replaced 
$$
\frac{2 \Delta}{N} I_n^{N,C} \Big( \frac{b_\mu \partial^2_\sigma a^2_\sigma}{a^4} - 2\frac{b_\mu (\partial_\sigma a^2_\sigma)^2}{a^6_\sigma} \Big) - \frac{\Delta^2}{N} \nu_n^{N,C} \Big( \frac{b^2_\mu \partial^2_\sigma a^2_\sigma}{a^4_\sigma} - 2\frac{b^2_\mu (\partial_\sigma a^2_\sigma)^2}{a^6_\sigma} \Big)
$$ 
with $o_{\mathbb{P}}(1)$ because of the convergences in probability gathered in Theorem \ref{th: convergence InN} and Proposition \ref{prop: converence nu}, respectively. It follows that $\frac{\Delta}{N} \partial^2_\sigma L_n^{N,C}(\theta)$}
converges in probability to 
$$\bar{\Pi} \Big( a^2_{\sigma_0} \Big( 2\frac{(\partial_\sigma a_\sigma^2)^2}{a_\sigma^6} - \frac{\partial^2_\sigma a_\sigma^2}{a^4_\sigma} \Big) \Big) + \bar{\Pi} \Big(\frac{\partial^2_\sigma a_\sigma^2}{a_\sigma^2}- \frac{(\partial_\sigma a_\sigma^2)^2}{a_\sigma^4} \Big) = : \Sigma^{(2)}(\theta)$$
uniformly in $\theta$, thanks to Proposition \ref{prop: converence nu} and Theorem \ref{th: convergence QnN}. We remark that $\Sigma^{(2)}(\theta) = \bar{\Pi}(\frac{(\partial_\sigma a^2_{\sigma})^2}{a_{\sigma}^4} )${\spa .} This implies that the matrix $\Sigma_n^{N,C}(\theta)$ converges in probability, uniformly in $\theta$, towards 
\begin{align*}
\Sigma (\theta) = \begin{pmatrix}
    2 \Sigma^{(1)}(\theta) & 0 \\
    0 &  \Sigma^{(2)} (\theta)
\end{pmatrix}.
\end{align*}
By the result above and \ref{as6} we know that the probability of $\int_0^1 \Sigma_n^{N,C}(\theta_0 + u ( \hat{\theta}_n^{N,C}- \theta_0)) du$ being invertible tends to $1$. Then, by applying its inverse to Equation \eqref{eq: start as norm} and using continuous mapping theorem, we conclude the proof for the complete observation case, remarking that the variance matrix is 
\begin{align*}
(\Sigma(\theta_0))^{-1}    
\times 
\begin{pmatrix}
4 \Sigma^{(1)}(\theta_0) & 0\\
0 & 2 \Sigma^{(2)}(\theta_0)  
\end{pmatrix} 
\times
(\Sigma(\theta_0))^{-1} = \begin{pmatrix}
( \Sigma^{(1)}(\theta_0) )^{-1} & 0 \\
0 & 2 (\Sigma^{(2)}(\theta_0))^{-1}
\end{pmatrix},
\end{align*} 
as we wanted.
\\
\\
Let us move to the partial observation framework. The proof follows the same lines as in the complete observation case with the only difference that, according to \eqref{eq: start partial 2}, we have 
$$ L_n^{N,P}(\theta) = \frac{3}{2} \frac{1}{\Delta} Q_n^{N,P} \Big( \frac{1}{a_\sigma^2} \Big) - 3  I_n^{N,P}\Big( \frac{b_\mu}{a_\sigma^2} \Big) + \frac{3}{2} \Delta \nu_n^{N,P} \Big( \frac{b_\mu}{a_\sigma^2}(b_\mu - 2 b_{\mu_0}) \Big) +  \nu_n^{N,P}\big(\log(a_\sigma^2)\big).$$
Then, Theorem \ref{th: convergence law InN} ensures that 
\begin{align*}
- \frac{1}{\sqrt{N}} \partial_\mu L_n^{N,P}(\theta_0) & = \frac{3}{\sqrt{N}}  I_n^{N,P}\Big(\frac{\partial_\mu b_{\mu_0}}{a^2_{\sigma_0}}\Big) \xrightarrow{d} \mathcal{N} \Big( 0, 9 \bar{\Pi} \Big( \Big(\frac{\partial_\mu b_{\mu_0}}{a_{\sigma_0}} \Big)^2 \Big) \Big).
\end{align*}
Moreover, Theorem \ref{th: convergence law QnN} implies 
$$- \sqrt{\frac{\Delta}{N} }\partial_\sigma L_n^{N,P}(\theta_0) \xrightarrow{d} \mathcal{N} \Big(0, \frac{9}{4} \bar{\Pi}\Big( \Big( \frac{\partial_\sigma a^2_{\sigma_0}}{a^2_{\sigma_0}} \Big)^2 \Big) \Big).$$
{\spa Recall that Theorem \ref{th: convergence QnN} applies with limiting factor $\frac{2}{3}$ in the partial observation case. Then,} acting as above, it is easy to check that, uniformly in $\theta$,
\begin{align*}
\Sigma_n^{N,P}(\theta) \xrightarrow{\mathbb{P}}
\begin{pmatrix}
     3 \Sigma^{(1)} (\theta) & 0\\
    0& \Sigma^{(2)}(\theta)
\end{pmatrix}
\end{align*}
with $\Sigma^{(j)} (\theta)$ defined as in the case of complete observation we dealt with before. 
We conclude again by applying the inverse of $\int_0^1 \Sigma_n^{N,P}(\theta_0 + u ( \hat{\theta}_n^{N,P}- \theta_0)) du$ in Equation \eqref{eq: start as norm} and using the continuous mapping theorem. The result is therefore proven. \\

\end{proof}

\appendix

\section{Proof {\spa of} technical results}{\label{s: proof technical}}

This section is dedicated to proving the lemmas and propositions introduced earlier. Their proofs are not straightforward; rather, they involve intricate technicalities. We will omit $n$ in the notation of the discretization step $\Delta_n \le 1$. We will write 
$b, a$ for $b_{\mu_0}, a_{\sigma_0}$, respectively. 


\subsection{Proof of Proposition \ref{p: error X tilde}} 
\begin{proof}
\noindent (i)
Observe that, by definition $\tilde X^i_{j\Delta} := \Delta^{-1} (Y^i_{(j+ 1)\Delta} - Y^i_{j\Delta})$ and that of $Y^i_s$, we have that
\begin{align*}
    \tilde X^i_{j\Delta} - X^i_{j\Delta} = \frac{1}{\Delta} \int_{j\Delta}^{(j+ 1)\Delta} (X_s^i - X^i_{j\Delta}) ds = B^i_j + A^i_j,
\end{align*}
where the dynamics of $X^i_s$ gathered in \eqref{eq:IPS} gives
\begin{align*}
B^i_j &:= \frac{1}{\Delta} \int_{j\Delta}^{(j+1)\Delta} \Big( \int_{j\Delta}^s b(Z^i_u, \Pi^N_u) du \Big) ds,\quad 
A^i_j := \frac{1}{\Delta}  \int_{j\Delta}^{(j+1)\Delta} \Big( \int_{j\Delta}^s a (Z^i_u, \Pi^N_u) d B^i_u \Big) d s.
\end{align*}
Application of the Fubini theorem then gives
$A^i_j = \alpha^i_j + \Delta^{\frac 1 2} a(Z^i_{j\Delta}, \Pi^N_{j\Delta}) \xi^i_j$ with $\xi^i_j$ defined in \eqref{eq: def xi tilde} and
\begin{align*}
\alpha^i_j &:= \frac{1}{\Delta} \int_{j\Delta}^{(j+1)\Delta} \Big( \int_{j\Delta}^s (a(Z^i_u, \Pi^N_u)-a(Z^i_{j\Delta},\Pi^N_{j\Delta})) d B^i_u \Big) ds\\ 
&= \frac{1}{\Delta} \int_{j\Delta}^{(j+1)\Delta} (a(Z^i_u, \Pi^N_u)-a(Z^i_{j\Delta},\Pi^N_{j\Delta})) ((j+1)\Delta - u) d B^i_u
\end{align*}
Let us show that 
$\varepsilon^i_j := \alpha^i_j + B^i_j$ 
has the stated property in (i): it suffices
to show that the collection of random variables $\Delta^{-1} \varepsilon^i_j$, where $0 < i \le N$, $0 \le j < n$, $n, N \in \mathbb{N}$, is bounded in $L^k$ for all $k \ge 2$. We decompose $\E [|\varepsilon^i_j|^k] \le C (\E [|\alpha^i_j|^k] + \E [|B^i_j|^k])$. 
Applying the 
Burkholder-Davis-Gundy and Jensen inequalities gives
\begin{align*}
\mathbb{E}
[|\alpha^i_j|^k] & \le \frac{C}{\Delta^k} \mathbb{E}
\Big[ \Big( \int_{j \Delta}^{(j+1)\Delta} ( a(Z^i_u,\Pi^N_u) - a(Z^i_{j\Delta},\Pi^N_{j\Delta}) )^2 ( (j+1) \Delta - u )^2 d u \Big)^{\frac k 2} \Big]\\
&\le \frac{C}{\Delta^{\frac{k}{2} + 1}} \int_{j\Delta}^{(j+1)\Delta} \mathbb{E}
[|a(Z^i_s, \Pi^N_s)-a(Z^i_{j\Delta},\Pi^N_{j\Delta})|^k ] ((j+1)\Delta-u)^k du,
\end{align*}
where
\begin{align*}
\mathbb{E}
[|a(Z^i_s, \Pi^N_s)-a(Z^i_{j\Delta},\Pi^N_{j\Delta})|^k ] \le C \mathbb{E}
[|Z^i_s-Z^i_{j\Delta}|^k+ W^k_2(\Pi^N_s,\Pi^N_{j\Delta})] \le C (s-j\Delta)^{\frac k 2}
\end{align*}
follows from the Lipschitz continuity of $a$ in \ref{as2} and Lemma \ref{l: moments}(ii), (iii).
We obtain $\mathbb{E} [|\alpha^i_j|^k] \le C \Delta^k$. 
Similarly, applying 
the Jensen inequality twice gives
\begin{align*}
\mathbb{E}
[|B^i_j|^k] &\le \Delta^{-1} \int_{j\Delta}^{(j+1)\Delta} \mathbb{E}
\Big[ \Big| \int_{j\Delta}^s b(Z^i_u, \Pi^N_u) d u \Big|^k \Big] d s\\ &\le \Delta^{k-2} \int_{j\Delta}^{(j+1)\Delta} \Big( \int_{j\Delta}^s \mathbb{E}
[ |b(Z^i_u, \Pi^N_u)|^k ] d u \Big) d s.
\end{align*}
The Lipschitz continuity of $b$ gathered in \ref{as2} implies
\begin{align}\label{ineq:condexp_b}
\mathbb{E}
[|b(Z^i_u, \Pi^N_u)|^k] &\le C  (|b(0, \delta_0)|^k + \mathbb{E}
[|Z^i_u|^k + W_2^k (\Pi^N_u, \delta_0)]),
\end{align}
where
$$
W^k_2 (\Pi^N_u, \delta_0) \le \Big( \frac{1}{N} \sum_{i=1}^N |Z^i_u|^2 \Big)^{\frac{k}{2}} \le \frac{1}{N} \sum_{i=1}^N |Z^i_u|^k
$$
by the Jensen inequality.
Then $\E [|b(Z^i_u,\Pi^N_u)|^k] \le C$ by Lemma \ref{l: moments}(i). We get
$\E
[|B^i_j|^k] \le
C \Delta^k
$,
which concludes the proof of part (i).\\
\\
(ii) It is enough to show that a collection of $\Delta^{-\frac{1}{2}}(\tilde X^i_{j\Delta} - X^i_{j\Delta})$, where $0 < i \le N$, $0 \le j < n$, $n,N \in \mathbb{N}$, is bounded in $L^k$ for all $k$.
According to (i) above, 
$$
|\tilde X^i_{j\Delta} - X^i_{j\Delta}|^k 
\le  C (\Delta^{\frac k 2} |a(Z^i_{j\Delta}, \Pi_{j \Delta}^N)|^k 
|\xi^i_j|^k
+ 
|\varepsilon^i_j|^k
), 
$$
where
$\E[|\varepsilon^i_j|^k] \le C \Delta^{k} \le C \Delta^{\frac{k}{2}}$
and
$\E [|a(Z^i_{j\Delta}, \Pi_{j \Delta}^N)|^k] \le C
$ 
similarly as in \eqref{ineq:condexp_b}, whereas
\begin{align*}
\mathbb{E} [|\xi^i_j|^k] = C \mathbb{E} [|\xi^i_j|^2]^{\frac{k}{2}} = C \Delta^{-\frac{3k}{2}} \Big( \int_{j\Delta}^{(j+ 1)\Delta}((j+1)\Delta - u)^2 du \Big)^{\frac k 2} = C
\end{align*}
for all $k$.
\\
\\
(iii) By the definition of Wasserstein distance, we have that
$$
W_2^k(\tilde{\Pi}_{j\Delta}^{N}, \Pi_{j\Delta}^{N}) \le \Big( \frac{1}{N} \sum_{i = 1}^N |\tilde Z^i_{j\Delta} - Z^i_{j\Delta}|^2 \Big)^{\frac k 2} = \Big( \frac{1}{N} \sum_{i = 1}^N |\tilde X^i_{j\Delta} - X^i_{j\Delta}|^2 \Big)^{\frac k 2},
$$
where applying the Jensen inequality gives
$$
W_2^k(\tilde{\Pi}_{j\Delta}^{N}, \Pi_{j\Delta}^{N}) \le \frac{1}{N} \sum_{i = 1}^N |\tilde X^i_{j\Delta} - X^i_{j\Delta}|^k.
$$
Hence, (ii) above ensures 
$\E [W_2^k(\tilde{\Pi}_{j\Delta}^{N}, \Pi_{j\Delta}^{N})] \le C \Delta^{\frac k 2}$,
which concludes the proof of our proposition.
\end{proof}

\subsection{Proof of Proposition \ref{prop: increments tilde}}

\begin{proof}
(i) By the definition of $\tilde X_{j\Delta}$ and that of $Y^i_s$, we have
\begin{align*}
\tilde X^i_{(j+1)\Delta} - \tilde X^i_{j\Delta} = \frac{1}{\Delta} \int^{(j+1)\Delta}_{j\Delta} (X^i_{s + \Delta} - X^i_s) d s = \tilde A^i_j + \tilde B^i_j,
\end{align*}
where we decompose the dynamics of $X^i_s$:
\begin{align*}
\tilde A^i_j &:= \frac{1}{\Delta} \int_{j\Delta}^{(j+ 1)\Delta} \Big( \int_{s}^{s + \Delta} a(Z_u^i, \Pi_u^N) dB_u^i \Big) ds,\\ 
\tilde B^i_j &:= \frac{1}{\Delta} \int_{j\Delta}^{(j+ 1)\Delta} \Big( \int_{s}^{s + \Delta} b(Z_u^i, \Pi_u^N) du \Big) ds.    
\end{align*}
We use the Fubini theorem on both terms, obtaining 
\begin{align*}
\tilde A_j^i & = \frac{1}{\Delta} \int_{j\Delta}^{(j+ 1)\Delta} a(Z_u^i, \Pi_u^N)(u - j\Delta) dB_u^i +  \frac{1}{\Delta} \int_{(j + 1)\Delta}^{(j+ 2)\Delta} a(Z_u^i, \Pi_u^N)((j + 2)\Delta - u) dB_u^i\\
&=  \Delta^{\frac{1}{2}} a(Z_{j\Delta}^i, \Pi_{j\Delta}^N) U_j^i + \tilde \alpha_j^i + \alpha_{j + 1}^i,
\end{align*}
where $U^i_j := \tilde \xi_j^i + \xi_{j + 1}^i$ as in \eqref{eq: def xi tilde}, \eqref{eq: def xi} and 
\begin{align*}
\tilde \alpha_j^i &:= \frac{1}{\Delta} \int_{j\Delta}^{(j+ 1)\Delta} (a(Z_u^i, \Pi_u^N) -  a(Z_{j\Delta}^i, \Pi_{j\Delta}^N))(u - j\Delta) dB_u^i,\\
\alpha_{j + 1}^i &:= \frac{1}{\Delta} \int_{(j + 1)\Delta}^{(j+ 2)\Delta} (a(Z_u^i, \Pi_u^N)-  a(Z_{j\Delta}^i, \Pi_{j\Delta}^N))((j + 2)\Delta - u) dB_u^i,
\end{align*}
and, similarly,
\begin{align*}
\tilde B_j^i 
=   \Delta b(\tilde{Z}_{j\Delta}^i, \tilde{\Pi}_{j\Delta}^N) + 
\tilde \beta_j^i + \beta_{j + 1}^i,
\end{align*}
where 
\begin{align*}
\tilde \beta_j^i &:= \frac{1}{\Delta} \int_{j\Delta}^{(j+ 1)\Delta} (b(Z_u^i, \Pi_u^N) - b(\tilde{Z}_{j\Delta}^i, \tilde{\Pi}_{j\Delta}^N))(u - j\Delta) du,\\
\beta_{j + 1}^i &:= \frac{1}{\Delta} \int_{(j + 1)\Delta}^{(j+ 2)\Delta} (b(Z_u^i, \Pi_u^N)- b(\tilde{Z}_{j\Delta}^i, \tilde{\Pi}_{j\Delta}^N))((j + 2)\Delta - u) du.
\end{align*}
With the notation $\tilde{\varepsilon}_j^i := \tilde \alpha_j^i + \alpha_{j + 1}^i + \tilde \beta_j^i + \beta_{j + 1}^i$, it implies 
\begin{align*}
\tilde X^i_{(j+ 1)\Delta} - \tilde X^i_{j\Delta} -  \Delta b(\tilde{Z}_{j\Delta}^i, \tilde{\Pi}_{j\Delta}^N) 
 =   \Delta^{\frac{1}{2}} a(Z_{j\Delta}^i, \Pi_{j\Delta}^N) U_j^i + \tilde{\varepsilon}_j^i.
\end{align*}
The proof of the first point is concluded once we show that $\tilde{\varepsilon}_j^i$ has the wanted properties. \\
Let us start by checking that $\E [|\E_{j\Delta}[ \tilde{\varepsilon}_j^i]|^k] \le C \Delta^{\frac{3}{2}k}$ for $k \ge 1$. Observe that $\E_{j\Delta}[ \tilde \alpha_j^i] = 0 = \E_{j\Delta}[ \alpha_{j + 1}^i]$. Hence, $|\E_{j\Delta}[ \tilde{\varepsilon}_j^i]|^k \le C \E_{j\Delta}[|\tilde \beta^i_j|^k +|\beta^i_{j+1}|^k ]$ by the Jensen inequality. We then use the Jensen inequality on $\tilde \beta_j^i$ 
obtaining that
\begin{align}{\label{eq: bound beta ij}}
\E
[ |\tilde \beta_j^i|^k] & \le \frac{C}{\Delta} \int_{j\Delta}^{(j+ 1)\Delta} (\E
[|b(Z_u^i, \Pi_u^N) - b({Z}_{j\Delta}^i, {\Pi}_{j\Delta}^N)|^k] \nonumber \\
& \qquad \qquad + \E
[|b({Z}_{j\Delta}^i, {\Pi}_{j\Delta}^N) - b(\tilde{Z}_{j\Delta}^i, \tilde{\Pi}_{j\Delta}^N)|^k]) (u - j\Delta)^k du \le C \Delta^{\frac{3}{2}k}
\end{align}
having also employed the Lipschitz continuity of $b$ and Lemma \ref{l: moments}(ii), (iii), 
Proposition \ref{p: error X tilde}(ii), (iii). Similarly, we obtain
\begin{align}{\label{eq: bound tilde beta i j+1}}
\E
[ |\beta_{j+ 1}^i|^k]  \le 
C \Delta^{\frac{3}{2}k}.
\end{align}
Let us now show that
$\E [ |\tilde{\varepsilon}_j^i|^k] \le C \Delta^k$ for $k \ge 2$. From the Burkholder-Davis-Gundy and Jensen inequalities we have 
\begin{align*}
\E
[ |\tilde \alpha_j^i|^k] & \le 
C \Delta^{-k + \frac{k}{2} - 1} \int_{j\Delta}^{(j+ 1)\Delta} \E
[|a(Z_u^i, \Pi_u^N) - a({Z}_{j\Delta}^i, {\Pi}_{j\Delta}^N)|^k] (u - j\Delta)^k du  \\
& \le C \Delta^{-k+\frac{k}{2}+\frac{k}{2}+k} = C \Delta^k
\end{align*}
having also used the Lipschitz continuity of $a$ and Lemma \ref{l: moments}(ii), (iii). 
A similar bound holds true for the moments of $\alpha_{j + 1}^i$ which, together with \eqref{eq: bound beta ij} and \eqref{eq: bound tilde beta i j+1}, provides 
$$\E
[ |\tilde{\varepsilon}_j^i|^k] \le C (\Delta^k + \Delta^{\frac{3k}{2}}) \le C \Delta^k.
$$
To conclude the proof of (i) it remains to deal with $\E [|\mathbb{E}_{j\Delta}[\tilde \varepsilon^i_j U^i_j]|^k ] \le 
\E [ | \tilde \varepsilon^i_j U^i_j |^k ]
\le \mathbb{E} [|\tilde \varepsilon^i_j|^{2k}]^{\frac 1 2} \mathbb{E} [| U^i_j|^{2k}]^{\frac 1 2}$, where
$\mathbb{E} [|\tilde \varepsilon^i_j|^{2k}] \le C \Delta^{2k}$ 
follows from above, whereas thanks to the independence and zero-mean normal distribution of the increments of the Brownian motion, 
\begin{equation}\label{eq: bound U 3.5}
\mathbb{E} [| U^i_j|^{2k}] = C \mathbb{E} [| U^i_j|^2]^k = C \E[|\tilde \xi^i_j|^2 + |\xi^i_{j+1}|^2]^k = C.
\end{equation}


Let us improve the rate of $|\E_{j\Delta} [\tilde \varepsilon^i_j U^i_j] |$. We consider the $k$-th absolute moment of 
$$
\mathbb{E}_{j\Delta} [ \tilde \varepsilon^i_j U^i_j ] = \mathbb{E}_{j\Delta} [ \tilde \alpha^i_j \tilde \xi^i_j + \tilde \beta^i_j \tilde \xi^i_j + \beta^i_{j+1} \tilde \xi^i_j + \beta^i_{j+1} \xi^i_{j+1}  + \alpha^i_{j+1} \xi^i_{j+1}].
$$ 
Firstly, we deal with that of $\mathbb{E}_{j\Delta} [ \tilde \beta^i_j \tilde \xi^i_j ]$. Let $k,l > 1$ be such that $\frac{1}{k}+\frac{1}{l}=1$. 
By H\"older's inequality,
$|\mathbb{E}_{j\Delta} [ \tilde \beta^i_j \tilde \xi^i_j ]| \le \mathbb{E}_{j\Delta} [|\tilde \beta^i_j|^k]^{\frac{1}{k}} \mathbb{E}_{j\Delta}[|\tilde \xi^i_j|^l]^{\frac{1}{l}}$,
where 
$\mathbb{E}_{j\Delta}[|\tilde \xi^i_j|^l] = \mathbb{E} [|\tilde \xi^i_j|^l] = C \mathbb{E} [|\tilde \xi^i_j|^2]^{\frac{l}{2}} = C$, we get that
$$
\mathbb{E} [|\mathbb{E}_{j\Delta} [ \tilde \beta^i_j \tilde \xi^i_j ]|^k ] \le C \mathbb{E} [|\tilde \beta^i_j|^k].
$$ 
By \eqref{eq: bound beta ij} we conclude that $\mathbb{E}[|\mathbb{E}_{j\Delta} [ \tilde \beta^i_j \tilde \xi^i_j ]|^k] \le C \Delta^{\frac{3}{2}k}$, where $C$ does not depend on $i,j,n,N$. The same argument works for the $k$-th absolute moment of $\mathbb{E}_{j\Delta}[\beta^i_{j+1} U^i_j]$. It remains to deal with
\begin{align*}
\mathbb{E}_{j\Delta} [\tilde \alpha^i_j \tilde \xi^i_j] &= \frac{1}{\Delta^{\frac{5}{2}}} \int_{j\Delta}^{(j+1)\Delta} \mathbb{E}_{j\Delta} [ (a(Z^i_s,\Pi^N_s)-a(\tilde Z^i_{j\Delta}, \tilde \Pi^N_{j\Delta})) ] (s-j\Delta)^2 d s,\\
\mathbb{E}_{j\Delta} [ \alpha^i_{j+1} \xi^i_{j+1} ] &= \frac{1}{\Delta^{\frac{5}{2}}} \int_{(j+1)\Delta}^{(j+2)\Delta} \mathbb{E}_{j\Delta} [(a(Z^i_s, \Pi^N_s)-a(\tilde Z^i_{j\Delta}, \tilde \Pi^N_{j\Delta}))] ((j+2)\Delta-s)^2 d s.
\end{align*}
Let us start considering the first. By Jensen's inequality it satisfies
$$
|\mathbb{E}_{j\Delta}[\tilde \alpha^i_j \tilde \xi^i_j]|^k \le \frac{1}{\Delta^{\frac{3}{2}k+1}} \int_{j\Delta}^{(j+1)\Delta} | \mathbb{E}_{j\Delta} [ a(Z^i_s,\Pi^N_s) - a(\tilde Z^i_{j\Delta},\tilde \Pi^N_{j\Delta}) ]|^k  |s-j\Delta|^{2k} ds,
$$
where
\begin{align*}
|\mathbb{E}_{j\Delta} [ a(Z^i_s,\Pi^N_s) - a(\tilde Z^i_{j\Delta},\tilde \Pi^N_{j\Delta}) ]|^k 
&\le 
C (|\mathbb{E}_{j\Delta} [ a(Z^i_s,\Pi^N_s)] - a(Z^i_{j\Delta},\Pi^N_{j\Delta}) |^p\\ 
&\qquad + |a(Z^i_{j\Delta},\Pi^N_{j\Delta}) - \mathbb{E}_{j\Delta} [ a(\tilde Z^i_{j\Delta},\tilde \Pi^N_{j\Delta}) ]|^p).
\end{align*}
Thanks to Proposition \ref{prop: ito}(ii) using the specific dependence form of $a$ on the probability measure variable, we get $\mathbb{E} [|\mathbb{E}_{j\Delta}[\tilde \alpha^i_j \tilde \xi^i_j]|^k] \le C \Delta^{\frac{3}{2}k}$ with the same $C$ for all $i,j,n,N$. The same argument works for the $k$-th absolute moment of $\mathbb{E}_{j\Delta}[\alpha^i_{j+1} \xi^i_{j+1}]$.\\

(ii)
By Proposition \ref{p: error X tilde}(ii), Lemma \ref{l: moments}(ii), we have that $\E [|\tilde X^i_{(j+ 1)\Delta} - \tilde X^i_{j\Delta}|^k] \le C \Delta^{\frac{k}{2}}$. 
\end{proof}

\subsection{Proof of Proposition \ref{prop: ito}}

\begin{proof}
Let us start considering Part (i). 
Write $z_i := (y_i,x_i) \in \R^2$ in the vector $\boldsymbol{z} := (z_1, \dots , z_N) \in \R^{2N}$ and note that under Assumption \ref{as7}, the function
$$\boldsymbol{z} \mapsto f \Big( z_1, \frac{1}{N} \sum_{i = 1}^N \delta_{z_i} \Big) = 
F \Big(z_1, \frac{1}{N}\sum_{i=1}^N K(z_1,z_i)\Big) =:
\phi(\boldsymbol{z}) 
$$
is twice continuously differentiable. 
Denote by $\partial^l_{x_i} \phi$ the $l$-th order partial derivative of $\phi$ with respect to $x_i$. Write 
$Z^i_s = (Y^i_s,X^i_s)$ in the vector $\boldsymbol{Z}_s := (Z^1_s, \dots, Z^N_s)$, where the processes $(Z_s^i)_{s \in [0, T]}$, $i=1,\dots,N,$ satisfy \eqref{eq:IPS}. Write $t_j := j \Delta_n$, $j=0,\dots, n$.

We now apply the multi-dimensional It\^o's formula 
as follows 
\begin{align*}
&f(Z^1_{t_j}, \Pi^N_{t_j}) - f(Z^1_{t_{j-1}}, \Pi^N_{t_{j-1}}) =
\phi (\boldsymbol{Z}_{t_j}) - \phi(\boldsymbol{Z}_{t_{j-1}}) \\
&\qquad= \sum_{i = 1}^N \Big( \int_{t_{j-1}}^{t_j} \partial_{y_i} \phi (\boldsymbol{Z}_s) dY_s^i + 
\int_{t_{j-1}}^{t_j}\partial_{x_i} \phi(\boldsymbol{Z}_s) dX_s^i + \frac{1}{2} \int_{t_{j-1}}^{t_j} \partial_{x_i}^{2} \phi(\boldsymbol{Z}_s) d \langle X^i \rangle_s \Big),
\end{align*}
where $\langle X^i \rangle$ denotes the quadratic variation of $(X^i_s)_{s \in [0,T]}$.
Furthermore, from the dynamics of the process $(X_s^i)_{t \in [0, T]}$ and $(Y_s^i)_{t \in [0, T]}$
this is 
\begin{align*}
&\sum_{i = 1}^N \Big( \int_{t_{j-1}}^{t_j}  \partial_{y_i} \phi (\boldsymbol{Z}_s) X_s^i ds + \int_{t_{j-1}}^{t_j} \partial_{x_i} \phi (\boldsymbol{Z}_s) b (Z^i_s, \Pi^N_s ) d s\\ 
&\qquad + \int_{t_{j-1}}^{t_j} \partial_{x_i} \phi (\boldsymbol{Z}_s) a (Z^i_s, \Pi^N_s) d B^i_s + \frac{1}{2} \int_{t_{j-1}}^{t_j} \partial_{x_i}^{2} \phi (\boldsymbol{Z}_s) a^2(Z^i_s, \Pi^N_s) ds\Big).
\end{align*}
Since the driving $(B_s^1, \dots , B_s^N)_{s \in [t_{j -1}, t_j]}$ is independent of $\mathcal{F}^N_{t_{j-1}}$, we obtain
\begin{equation}{\label{eq: ito 1}}
\E_{t_{j - 1}} [f (Z^1_{t_j}, \Pi_{t_j}^N)] -  f (Z^1_{t_{j-1}}, \Pi_{t_{j-1}}^N) = \E_{t_{j - 1}} \Big[ \int_{t_{j-1}}^{t_j} A \phi(\boldsymbol{Z}_s) d s \Big],
\end{equation}
where
$$
A \phi (\boldsymbol{Z}_s) := \sum_{i = 1}^N ( \partial_{y_i} \phi (\boldsymbol{Z}_s) X_s^i + \partial_{x_i} \phi (\boldsymbol{Z}_s) b (Z^i_s, \Pi^N_s ) + \frac{1}{2} \partial_{x_i}^2 \phi (\boldsymbol{Z}_s) a^2(Z^i_s, \Pi^N_s )). 
$$
To conclude, we need to bound $\partial_{y_i} \phi$, $\partial_{x_i} \phi$, $\partial_{x_i}^{2} \phi$, that are actually the derivatives of $f$.
To do that, we rely on the assumption about the dependence of $f$ on the convolution with a probability measure. To compute the derivative with respect to $x_i$, 
we have to consider two different cases, depending on whether $i \neq 1$ or $i = 1$. For $i \neq 1$,
$$
\partial_{x_i} \phi (\boldsymbol{z}) = 
\partial_w F \Big( z_1, \frac{1}{N} \sum_{l = 1}^N K(z_1, z_l) \Big) \frac{1}{N} \partial_{x_2} K(z_1, z_i),
$$
where $\partial_w F (z,w)$ denotes the partial derivative of a function $F(z,w)$ of two variables $(z,w) \in \R^2 \times \R$ with respect to the second one and $\partial_{x_2} K (z_1,z_2)$ denotes the partial derivative of a function $K(z_1,z_2)$ with respect to the component $x_2$ of its second variable $z_2 = (y_2,x_2) \in \R^2$. 
For $i=1$, instead, denoting as $\partial_x F (z,w)$ the partial derivative of $F(z,w)$ with respect to the component $x$ of the variable $z := (y,x)$,
we have 
\begin{align*}
&\partial_{x_1} \phi (\boldsymbol{z})  = \partial_{x} F \Big( z_1, \frac{1}{N} \sum_{l = 1}^N K(z_1, z_l) \Big) \\
&\qquad + \partial_w F \Big( z_1, \frac{1}{N} \sum_{l = 1}^N K(z_1, z_l) \Big) \frac{1}{N} \sum_{l = 1}^N ( \partial_{x_1} K (z_1, z_l) + \partial_{x_2} K (z_1, z_l) ).
\end{align*}
An analogous result hold for $\partial_{y_i} \phi$, $\partial_{x_i}^2 \phi$. 
We recall that all the partial derivatives of $F$, $K$ up to the second order have polynomial growth. 
Then we can see that
$A \phi(\boldsymbol{Z}_s)$ appearing in \eqref{eq: ito 1} is bounded in $L^p$, $p \ge 1$, uniformly in $s \in [0,T]$, $N \in \mathbb{N}$. Hence, by Jensen's inequality for conditional expectation and Minkowski's inequality, 
\begin{align*}
&\E [ | \E_{t_{j-1}}[f(Z^1_{t_j},\Pi^N_{t_j})] - f(Z^1_{t_{j-1}},\Pi^N_{t_{j-1}}) |^p ]^{\frac{1}{p}} \\ 
&\qquad  \le \E \Big[ \Big| \int_{t_{j-1}}^{t_j} A \phi (\boldsymbol{Z}_s) ds \Big|^p \Big]^{\frac{1}{p}} \le \int_{t_{j-1}}^{t_j} \E [|A\phi(\boldsymbol{Z}_s)|^p]^{\frac{1}{p}} ds \le C \Delta_n,
\end{align*}
as we wanted. 

Let us move to the proof of Part (ii). Taylor's theorem gives the approximation:
\begin{align*}
&f(\tilde Z^1_{t_{j-1}},\tilde \Pi^N_{t_{j-1}}) - f(Z^1_{t_{j-1}},\Pi^N_{t_{j-1}}) = \phi (\boldsymbol{\tilde Z}_{t_{j-1}}) - \phi (\boldsymbol{Z}_{t_{j-1}})\\ 
&\qquad= \sum_{k=1}^N \partial_{x_k} \phi (\boldsymbol{Z}_{t_{j-1}}) (\tilde X^k_{t_{j-1}} - X^k_{t_{j-1}})  + \sum_{k,l=1}^N R_{k l} (\boldsymbol{\tilde Z}_{t_{j-1}}) (\tilde X^k_{t_{j-1}}- X^k_{t_{j-1}}) (\tilde X^l_{t_{j-1}}- X^l_{t_{j-1}}),
\end{align*}
where 
$$
R_{kl} (\boldsymbol{\tilde Z}_{t_{j-1}}) := \int_0^1 (1-s) \partial_{x_k} \partial_{x_l} \phi (\boldsymbol{Z}_{t_j} + s (\boldsymbol{\tilde Z}_{t_{j-1}}-\boldsymbol{Z}_{t_{j-1}})) d s
$$
for all $k,l = 1, \dots, N$. Let $p \ge 1$. By Minkowski's inequality, we have that
\begin{align*}
&\E [ | \E_{t_{j-1}} [f(\tilde Z^1_{t_{j-1}}, \tilde \Pi^N_{t_{j-1}})] - f(Z^1_{t_{j-1}}, \Pi^N_{t_{j-1}}) |^p ]^{\frac{1}{p}}\\ 
&\qquad \le \sum_{k=1}^N \E [ | \partial_{x_k} \phi (\boldsymbol{Z}_{t_{j-1}}) (\E_{t_{j-1}} [\tilde X^k_{t_{j-1}}] - X^k_{t_{j-1}}) |^p ]^{\frac{1}{p}}\\
&\qquad \qquad + \sum_{k,l=1}^N \E [| R_{kl} (\boldsymbol{\tilde Z}_{t_{j-1}}) (\tilde X^k_{t_{j-1}}-X^k_{t_{j-1}}) (\tilde X^l_{t_{j-1}} - X^l_{t_{j-1}}) |^p ]^{\frac{1}{p}},
\end{align*}
where Proposition \ref{p: error X tilde}(ii) implies 
$\E [| \tilde X^k_{t_{j-1}} - X^k_{t_{j-1}} |^p ]^{\frac{1}{p}} \le C_p \Delta_n^{\frac{1}{2}}$, whereas $\E [|\E_{t_{j-1}} [\tilde X^k_{t_{j-1}}] - X^k_{t_{j-1}}|^p]^{\frac{1}{p}} \le C_p \Delta_n$ follows from
$$
\tilde X^k_{t_{j-1}} - X^k_{t_{j-1}} = \frac{1}{\Delta} \int_{t_{j-1}}^{t_j}  (X^k_s-X^k_{t_{j-1}}) ds
$$
using the dynamics of $(X^k_s)_{s \in [0,T]}$ in the very same way as in Part (i). Also, it is easy to check that $\sum_{k=1}^N \E[ |\partial_{x_k} \phi (\boldsymbol{Z}_{t_{j-1}})|^p]^{\frac{1}{p}} \le C_p$ and $\sum_{k,l=1}^N \E [|R_{k,l}(\boldsymbol{\tilde Z_{t_{j-1}}})|^p]^{\frac{1}{p}} \le C_p$. Since $p$ is arbitrary, we obtain
$$
\E [ | \E_{t_{j-1}} [f(\tilde Z^1_{t_{j-1}}, \tilde \Pi^N_{t_{j-1}})] - f(Z^1_{t_{j-1}}, \Pi^N_{t_{j-1}}) |^p ]^{\frac{1}{p}} \le C_p \Delta_n,
$$
which concludes the proof.

\end{proof}

\subsection{Proof of Proposition \ref{prop: converence nu}}
\begin{proof}
For every $\theta$, the convergence $\frac{\Delta}{N} \nu_n^{N,C}(f_\theta) \xrightarrow{\mathbb{P}} \bar{\Pi}(f_\theta)$ follows from Lemma~5.2 of \cite{McKean}, which uses the condition \eqref{eq: cond pol growth}. The tightness in the space of continuous functions on $\Theta$ follows in the same way as that of $\theta \mapsto I^N_n(\theta)$ defined in (18) in the proof of (15) in \cite{McKean} by using the polynomial growth of $\sup_{\theta} |\nabla_\theta f_\theta|$ and uniform in $i, j, n, N$ bound on the moments of $Z_{j\Delta}^i$ established in Lemma~\ref{l: moments}(i). 
It remains to prove that the same convergence holds true if we replace $Z^i_{j\Delta}$ with $\tilde Z^i_{j\Delta}$. In the same way we can see the tightness, whereas the pointwise convergence follows if we prove that for every $\theta$,
\begin{equation*}\label{eq: nuPdecomp}
    \frac{\Delta}{N}(\nu_n^{N,P}(f_\theta) - \nu^{N,C}_n (f_\theta)) = \frac{\Delta}{N} \sum_{i=1}^N \sum_{j=0}^{n-1} (f_\theta (\tilde Z^i_{j\Delta}, \tilde \Pi^N_{j\Delta}) - f_\theta (Z^i_{j\Delta}, \Pi^N_{j\Delta}))
\end{equation*}
converges to $0$ in $L^1$. By \eqref{eq: cond pol growth}, 
\begin{align}{\label{eq: increments f}}
    &\E [ |f_\theta (\tilde Z^i_{j\Delta}, \tilde \Pi^N_{j\Delta}) - f_\theta (Z^i_{j\Delta}, \Pi^N_{j\Delta})| ] \le C \E \big[ (|\tilde Z^i_{j\Delta}-Z^i_{j\Delta}| + W_2(\tilde \Pi^N_{j\Delta},\Pi^N_{j\Delta}))\\
    &\qquad \times (1+|\tilde Z^i_{j\Delta}|^k + |Z^i_{j\Delta}|^k+W_2^k(\tilde \Pi^N_{j\Delta},\delta_0)+W_2^k(\Pi^N_{j\Delta},\delta_0)) \big]. 
    \nonumber
\end{align}
The Cauchy-Schwarz inequality applies to the product on the RHS of \eqref{eq: increments f}. Consider the $L^2$-norm of its factors. The first one satisfies 
$$
\E \big[ |\tilde Z^i_{j\Delta}-Z^i_{j\Delta}|^2 + W_2^2 (\tilde \Pi^N_{j\Delta},\Pi^N_{j\Delta}) \big] \le C \Delta
$$ for all $i,j,n,N$ by Proposition \ref{p: error X tilde}(ii), (iii). 
Let us turn to the second factor on the RHS of \eqref{eq: increments f}. We note that $W_2^{2k} (\tilde \Pi^N_{j\Delta}, \delta_0) \le \frac{1}{N} \sum_{i=1}^N |\tilde Z^i_{j\Delta}|^{2k}$, whereas $\E [|\tilde Z^i_{j\Delta}|^{2k}] \le C$ for all $i,j,n,N$ follows from the uniform bound on the moments of $\tilde Z^i_{j\Delta}-Z^i_{j\Delta}$ and $Z^i_{j\Delta}$ by Proposition~\ref{p: error X tilde}(ii) and Lemma \ref{l: moments}(i) respectively. 
Hence,
$$
\E [ |
\frac{\Delta}{N}  (\nu_n^{N,P}(f_\theta) - \nu^{N,C}_n (f_\theta)) | ] \le C \Delta^{\frac{1}{2}}
$$
goes to $0$, completing the proof. 
\end{proof}

\subsection{Proof of Theorem \ref{th: convergence InN}}
\begin{proof} 
For a function $f : \mathbb{R}^2 \times \mathcal{P}_2(\mathbb{R}^2) \to \mathbb{R}$, we denote $f (Z^i_t, \Pi^N_t)$ and $f (\tilde Z^i_t, \tilde \Pi^N_t)$ by $f^i_t$ and $\tilde f^i_t$ respectively.
We use the same approach to deal with functionals $I_n^{N,C}(f_\theta)$ and $I_n^{N,P}(f_\theta)$. Since the proof for $I_n^{N,C}(f_\theta)$ is simpler, we will only present the proof for $I_n^{N,P}(f_\theta)$ here.
We have that
$$I_n^{N,P}(f_\theta) = \tilde{I}_n^{N,P}(f_\theta) + \Delta \sum_{j = 1}^{n -2} \sum_{i = 1}^N (\tilde f_\theta)_{(j- 1)\Delta}^i (\tilde b_{j\Delta}^i - \tilde b_{(j-1)\Delta}^i),
$$
where
\begin{equation}{\label{eq: def I tilde}}
\tilde{I}_n^{N, P }(f_\theta) := \sum_{j = 1}^{n -2} \sum_{i = 1}^N 
(\tilde{f}_\theta)_{(j- 1)\Delta}^i
(\tilde{X}_{(j+1)\Delta}^i - \tilde{X}_{j\Delta}^i - \Delta \tilde b_{j\Delta}^i).
\end{equation}
It suffices to prove that uniformly in $\theta$,
\begin{equation}\label{lim:tildeI}
\frac{1}{N} \tilde I^{N,P}_n (f_\theta) \xrightarrow{\mathbb{P}} 0, \quad \frac{1}{N} (I^{N,P}_n (f_\theta) -  \tilde I^{N,P}_n (f_\theta)) \xrightarrow{\mathbb{P}} 0.   
\end{equation}

We start from the LHS of \eqref{lim:tildeI}. We follow Proposition \ref{prop: increments tilde}(i) by writing $\tilde{X}_{(j+1)\Delta}^i - \tilde{X}_{j\Delta}^i - \Delta \tilde b_{j\Delta}^i=\Delta^{\frac{1}{2}} a^i_{j\Delta} U^i_j + \tilde{\varepsilon}_j^i$ in $\tilde I^{N,P}_n (f_\theta)$. Then it suffices to prove that for $k=1,2$, uniformly in $\theta$, 
\begin{equation}\label{lim:sumFk}
S_k (\theta) := \sum_{j=1}^{n-2} F_{j,k} (\theta) \xrightarrow{\mathbb{P}} 0,
\end{equation}
where 
\begin{equation}\label{def:Fk}
F_{j,1}(\theta) := \frac{\Delta^{\frac{1}{2}} }{N} \sum_{i = 1}^N (\tilde f_\theta)^i_{(j- 1)\Delta} a^i_{j\Delta} U^i_j,\quad F_{j,2}(\theta) := \frac{1}{N} \sum_{i = 1}^N (\tilde f_\theta)_{(j- 1)\Delta}^i \tilde{\varepsilon}_j^i.   
\end{equation}

In order to prove \eqref{lim:sumFk} for every $\theta$ we use \cite[Lemma 9]{GenJac93}. We note that $F_{j,k} (\theta)$ is $\mathcal{F}^N_{(j+2)\Delta}$-measurable, since so is $\tilde{X}_{(j+1)\Delta}^i$. 
Therefore, we apply \cite[Lemma 9]{GenJac93}  separately to the sum of $F_{j,k}(\theta)$ over even $j$ and that one over odd $j$. However, if we gather all the terms back, then it suffices to prove that, for $k=1,2$, 
\begin{align}{\label{eq: conv zjk n+2}}
\sum_{j = 1}^{n - 2} \mathbb{E}_{j\Delta}[F_{j,k}(\theta)] \xrightarrow{\mathbb{P}} 0,\quad
\sum_{j = 1}^{n - 2} \mathbb{E}_{j\Delta}[(F_{j,k}(\theta))^2] \xrightarrow{\mathbb{P}} 0. 
\end{align}

Denote the sum on the LHS of \eqref{eq: conv zjk n+2} by $\Sigma_k$ for $k=1,2$. Actually, we have $\Sigma_1 = 0$ since $\mathbb{E}_{j \Delta}[U_j^i] = 0$ 
and  
$(\tilde f_\theta)^i_{(j- 1)\Delta} a^i_{j\Delta}$ is $\mathcal{F}^N_{j\Delta}$-measurable.
Consider $\Sigma_2$. Polynomial growth of $f_\theta$, bounded moments of $\tilde Z^i_{(j-1)\Delta}$ and Proposition \ref{prop: increments tilde}(i) give 
$
\mathbb{E}[| (\tilde f_\theta)^i_{(j-1)\Delta} \mathbb{E}_{j \Delta} [\tilde{\varepsilon}_j^i] |] \le C \Delta^{\frac 3 2}
$ 
for all $i,j,n,N$.
Hence, $\mathbb{E}[|\Sigma_2|] \le C \Delta^{\frac 1 2}$, which proves the convergence on the LHS of \eqref{eq: conv zjk n+2} for $k=2$.

Let us turn to the RHS of \eqref{eq: conv zjk n+2}. For $k=1$, we get
\begin{align}\label{eq: comp square n+3,5}
\sum_{j = 1}^{n - 2} \mathbb{E}_{j\Delta}[(F_{j,1}(\theta))^2] = C \frac{\Delta}{N^2} \sum_{j = 1}^{n - 2} \sum_{i = 1}^N
(\tilde f_\theta^2)^i_{(j-1)\Delta} (a^2)^i_{j\Delta}, 
\end{align}
because
$\mathbb{E}_{j\Delta} [U^{i_1}_j U^{i_2}_j] = \mathbb{E} [U^{i_1}_j U^{i_2}_j] = C \mathbf{1} (i_1 \neq i_2)$
using independence, $\mathbb{E} [U^{i_1}_j] = \mathbb{E} [U^{i_2}_j] = 0$ and \eqref{eq: bound U 3.5}. Polynomial growth of both $f_\theta$, $a$ and bounded moments of $Z^i_{j\Delta}$ and, hence, those of $\tilde Z^i_{(j-1)\Delta}$ 
give the upper bound $C/N \to 0$ on the $L^1$-norm of the sum in \eqref{eq: comp square n+3,5}, which proves its convergence to $0$ in probability. 
For $k=2$, we get 
\begin{align}\label{eq: sumZ22}
\sum_{j = 1}^{n - 2} \mathbb{E}_{j\Delta}[(F_{j,2}(\theta))^2] &= \frac{1}{N^2} \sum_{j = 1}^{n - 2} \sum_{i_1, i_2 = 1}^N (\tilde f_\theta)_{(j- 1)\Delta}^{i_1}  (\tilde f_\theta)_{(j- 1)\Delta}^{i_2}
\mathbb{E}_{j\Delta}[\tilde{\varepsilon}^{i_1}_j \tilde{\varepsilon}^{i_2}_j],
\end{align}
where $\mathbb{E} [|\mathbb{E}_{j\Delta}[\tilde{\varepsilon}^{i_1}_j \tilde{\varepsilon}^{i_2}_j]|^k] \le \mathbb{E} [|\tilde{\varepsilon}^{i_1}_j \tilde{\varepsilon}^{i_2}_j |^k]$, moreover, $\mathbb{E} [(\tilde{\varepsilon}^{i}_j)^{2k}]^{\frac 1 2} \le C \Delta^k$ for $k\ge 1$ by Proposition~\ref{prop: increments tilde}(i).
Hence, the $L^1$-norm of the sum in \eqref{eq: sumZ22} is bounded by $C\Delta$.  
It completes the proof that \eqref{eq: sumZ22} goes to $0$ in probability
for every $\theta$. 

Concerning the uniformity in $\theta$, we need to prove the tightness of $\theta \mapsto S_k (\theta)$ in the space of continuous functions on $\Theta$ for $k=1,2$.
Let us start from $k=2$, as it is simpler to deal with. By Theorem 14.5 in \cite{40McK}, it suffices to prove that, for all $n,N$, 
\begin{equation}{\label{eq: unif z2 n+4}}
\mathbb{E} \Big[ \sup_{\theta} | \nabla_\theta S_2 (\theta) | \Big] \le C.
\end{equation}
We have that
$$\nabla_\theta  
S_2 (\theta) = \frac{1}{N} \sum_{i = 1}^N \sum_{j = 1}^{n - 2} (\nabla_\theta \tilde f_\theta)^i_{(j-1)\Delta} \tilde{\varepsilon}_j^i,$$ 
where $\sup_\theta |\nabla_\theta f_\theta|$ has a polynomial growth and the moments of  $\tilde Z^i_{(j-1)\Delta}$ and $\Delta^{-1} \tilde{\varepsilon}_j^i$ are uniformly bounded in $i,j,n,N$. It clearly implies \eqref{eq: unif z2 n+4}. 

Let us move to the uniformity for $k=1$. By Theorem 20 in Appendix 1 of \cite{39McK}, it suffices to prove that, for all $n,N$ and for all $\theta, \theta' \in \Theta$, 
\begin{equation}\label{ineq:tightsumF1}
\mathbb{E} [ | S_1(\theta) |^2 ] \le C, \qquad \mathbb{E} [ | S_1(\theta)-S_1(\theta') |^2 ] \le C| \theta - \theta'|^2.
\end{equation} 
Let us start from the LHS of \eqref{ineq:tightsumF1}. We note that for all $i_1 = i_2$ and $j_1 +2 \le j_2$, 
\begin{equation}\label{eq: mixed terms n+5}
\mathbb{E}[(\tilde f_\theta)_{(j_1-1)\Delta}^{i_1} a^{i_1}_{j_1\Delta} U^{i_1}_{j_1} (\tilde f_\theta)_{(j_2- 1)\Delta}^{i_2} a^{i_2}_{j_2\Delta}  U^{i_2}_{j_2}] =0
\end{equation}
follows from $\mathbb{E}_{j_2 \Delta}[U^{i_2}_{j_2}] = 0$, whereas independence of Brownian motions implies \eqref{eq: mixed terms n+5} for all $i_1 \neq i_2$ and $j_1, j_2$. We get that 
\begin{align}
\mathbb{E} [ | S_1 (\theta) |^2 ]  &= \frac{\Delta}{N^2} \sum_{i= 1}^N \sum_{j_1, j_2 = 1}^{n - 2} \mathbf{1}( |j_1-j_2|<2)\nonumber\\ 
&\qquad \times \mathbb{E} [(\tilde f_\theta)_{(j_1-1)\Delta}^i a^i_{j_1\Delta} U^i_{j_1} (\tilde f_\theta)_{(j_2- 1)\Delta}^i a^i_{j_2\Delta}  U^i_{j_2}].\label{ineq: Z1diff2}
\end{align}
We note the polynomial growth of both $a$ and $\sup_\theta |f_\theta|$, because the mean value theorem gives $\sup_\theta |f_\theta| \le |\theta-\theta'| \sup_\theta |\nabla_\theta f_\theta| + |f_{\theta'}|$ and $\Theta$ is bounded. Moreover, the moments of $\tilde Z^i_{(j-1)\Delta}$ and $Z^i_{j\Delta}$ are bounded uniformly in $i,j,n,N$, whereas $\mathbb{E}[| U^i_j |^2] = C$ thanks to \eqref{eq: bound U 3.5}. Hence, we conclude that \eqref{ineq: Z1diff2} is bounded by $C/N$. Let us turn to the RHS of \eqref{ineq:tightsumF1}. It is straightforward that one can work in the same way replacing the function $f_\theta$ here above with the difference $f_\theta- f_{\theta'}$ and using the mean value theorem. Hence, the uniformity in $\theta$ follows.

It remains to prove the convergence on the RHS of \eqref{lim:tildeI}. It is enough to prove that
\begin{equation}\label{lim:I-tildeI}
\mathbb{E} \Big[ \sup_\theta \Big| \frac{1}{N} (I^{N,P}_n(f_\theta) - \tilde I^{N,P}_n(f_\theta)) \Big| \Big] \to 0.
\end{equation}
We recall that $\sup_\theta |f_\theta|$ has a polynomial growth, whereas $b$ is Lipschitz continuous according to \ref{as2}. Using the Cauchy-Schwarz inequality and Proposition \ref{prop: increments tilde}(ii),  
the LHS of \eqref{lim:I-tildeI} gets bounded by $C\Delta^{\frac 1 2}$.
It concludes the proof of the theorem.
\end{proof}

\subsection{Proof of Theorem \ref{th: convergence QnN}}
\begin{proof}
For a function $f : \mathbb{R}^2 \times \mathcal{P}_2(\mathbb{R}^2) \to \mathbb{R}$, we denote $f (Z^i_t, \Pi^N_t)$ and $f (\tilde Z^i_t, \tilde \Pi^N_t)$ by $f^i_t$ and $\tilde f^i_t$ respectively.
As in the proof of previous theorem, we decide to detail only the partial observation case. We aim at proving that uniformly in $\theta$,
\begin{equation}{\label{eq: start Qn n+7}}
\frac{1}{N} Q_n^{N,P}(f_\theta)
\xrightarrow{\mathbb{P}} \frac{2}{3} \bar{\Pi}(f_\theta a^2). \nonumber
\end{equation}
According to Proposition \ref{prop: increments tilde}(i) we have $\tilde{X}_{(j + 1)\Delta}^i - \tilde{X}_{j\Delta}^i = \Delta^\frac{1}{2} a_{j\Delta}^i U^i_j + \Delta \tilde b_{j\Delta}^i + \tilde \varepsilon^i_j$. We use it to decompose
$$
\frac{1}{N} Q_n^{N,P}(f_\theta) = \sum_{j=1}^{n-2} (\rho_{j,0} (\theta) + \rho_{j,1} (\theta)),
$$
where
\begin{align*}
    \rho_{j,1} (\theta) &:= \frac{\Delta}{N} \sum_{i=1}^N ( \tilde f_\theta)_{(j-1)\Delta}^i (a^2)_{j\Delta}^i (U^i_j)^2,\\
    \rho_{j,0} (\theta) &:= \frac{1}{N} \sum_{i=1}^N ( \tilde f_\theta)_{(j-1)\Delta}^i ( 2 \Delta^{\frac{1}{2}} a_{j\Delta}^i U^i_j + \Delta \tilde b_{j\Delta}^i + \tilde \varepsilon^i_j ) (\Delta \tilde b_{j\Delta}^i + \tilde \varepsilon^i_j).
\end{align*}
We are going to prove that uniformly in $\theta$,
\begin{equation}\label{lim:omega}
\sum_{j=1}^{n-2} \rho_{j,1} (\theta) \xrightarrow{\mathbb{P}} \frac{2}{3} \bar{\Pi}(f_\theta a^2), \qquad 
\sum_{j=1}^{n-2} \rho_{j,0} (\theta) \xrightarrow{\mathbb{P}} 0.
\end{equation}

In order to prove the convergence on the LHS of \eqref{lim:omega} for every $\theta$, we use \cite[Lemma 9]{GenJac93} in the same way as in the proof of \eqref{lim:sumFk}: it suffices to prove that
\begin{equation}\label{lim:omega1}
\sum_{j = 1}^{n -2} \mathbb{E}_{j\Delta}[\rho_{j,1}(\theta)] \xrightarrow{\mathbb{P}} \frac{2}{3} \bar{\Pi}(f_\theta a^2), \qquad \sum_{j = 1}^{n -2} \mathbb{E}_{j\Delta}[(\rho_{j,1}(\theta))^2] \xrightarrow{\mathbb{P}} 0.    
\end{equation}
We have 
\begin{align*}
\sum_{j = 1}^{n -2} \mathbb{E}_{j\Delta}[\rho_{j,1}(\theta)] & = \frac{\Delta}{N} \sum_{j = 1}^{n -2} \sum_{i = 1}^N (\tilde f_\theta)_{(j- 1)\Delta}^i (a^2)_{j \Delta}^i \mathbb{E}_{j\Delta}[(U^i_j)^2],
\end{align*}
where $\mathbb{E}_{j\Delta}[(U^i_j)^2] = \frac{2}{3}$. Moreover, $a^2$ has the polynomial growth, moments of $Z^i_{j\Delta}$ are bounded, whereas, for every $p \ge 1$, the approximation
\begin{align*}
\mathbb{E} [|(\tilde f_\theta)_{(j- 1)\Delta}^i- (f_\theta)^i_{j \Delta}|^p]^{\frac{1}{p}} &\le \mathbb{E} [|(\tilde f_\theta)_{(j- 1)\Delta}^i- (f_\theta)^i_{(j-1) \Delta}|^p]^{\frac{1}{p}}\\ 
&\qquad+ \mathbb{E} [|(f_\theta)_{(j- 1)\Delta}^i- (f_\theta)^i_{j \Delta}|^p]^{\frac{1}{p}} \le C \Delta^{\frac{1}{2}}  
\end{align*}
holds true by \eqref{eq: cond pol growth} and Proposition \ref{prop: increments tilde}(ii). In the same way as in the proof of Proposition \ref{prop: converence nu} for complete observations, we get
$$
\frac{2}{3} \frac{\Delta}{N} \sum_{j=1}^{n-2} \sum_{i=1}^N (f_\theta)^i_{j\Delta} (a^2)^i_{j\Delta} \xrightarrow{\mathbb{P}} \frac{2}{3} \bar \Pi (f_\theta a^2),
$$
which completes the proof of the first convergence in \eqref{lim:omega1}. The second one in \eqref{lim:omega1} follows from the bound $C\Delta$ on the $L^1$-norm of
\begin{align*}
\sum_{j = 1}^{n -2} \mathbb{E}_{j\Delta}[(\rho_{j,1}(\theta))^2]  &= \Big( \frac{\Delta}{N} \Big)^2 \sum_{j = 1}^{n -2} \sum_{i_1, i_2 = 1}^N (\tilde f_\theta)_{(j- 1)\Delta}^{i_1} (a^2)_{j\Delta}^{i_1} (\tilde f_\theta) _{(j-1)\Delta}^{i_2} (a^2)_{j\Delta}^{i_2}\\ 
&\qquad \qquad \qquad \times  \mathbb{E}_{j\Delta}[(U^{i_1}_j)^2(U^{i_2}_j)^2],
\end{align*}
that goes to $0$ for $n \rightarrow \infty$. This finishes the proof of the convergence on the LHS of \eqref{lim:omega} for every $\theta$.

Let us prove that the convergence on the RHS of \eqref{lim:omega} in $L^1$ for every $\theta$. Indeed, having used Proposition \ref{prop: increments tilde}(i) and the polynomial growth of $f_\theta$, $b$ and $a$, we get that for $p \ge 1$,
$$
\mathbb{E}[|a^i_{j\Delta} U^i_j|^p]^{\frac{1}{p}} \le C, \qquad \mathbb{E} [|\Delta \tilde b^i_{j\Delta} + \tilde \varepsilon^i_j|^p]^{\frac{1}{p}} \le C \Delta, 
$$
which implies 
\begin{align}
\mathbb{E} \Big[ \Big| \sum_{j = 1}^{n -2} \rho_{j,0}(\theta) \Big| \Big]
&\le C \Delta^{\frac{1}{2}}.
\end{align}

The uniformity follows from 
\begin{align*}
&\mathbb{E} \Big[ \sup_\theta \Big|\frac{1}{N} \nabla_\theta Q_n^{N,P} (f_\theta) \Big| \Big]\\  
&\qquad \le  \frac{1}{N} \sum_{j=1}^{n-2} \sum_{i=1}^N
\mathbb{E} [ \sup_\theta |(\nabla_\theta \tilde f_\theta)^i_{(j-1)\Delta} | (\tilde X^i_{(j+1)\Delta} - \tilde X^i_{j\Delta})^2 ]\\
&\qquad \le \frac{1}{N} \sum_{j=1}^{n-2} \sum_{i=1}^N
\mathbb{E}[\sup_\theta |(\nabla_\theta \tilde f_\theta)^i_{(j-1)\Delta} |^2]^{\frac{1}{2}} \mathbb{E} [(\tilde X^i_{(j+1)\Delta} - \tilde X^i_{j\Delta})^4]^{\frac{1}{2}} \le C,
\end{align*}
where we have used the Cauchy-Schwarz inequality, polynomial growth of $\sup_\theta |\nabla_\theta f_\theta|$, bounded moments of $\tilde Z^i_{(j-1)\Delta}$ and Proposition \ref{prop: increments tilde}(ii). The proof is therefore concluded.  
\end{proof}

\subsection{Proof of Theorem \ref{th: convergence law InN}}
\begin{proof}
For a function $f : \mathbb{R}^2 \times \mathcal{P}_2(\mathbb{R}^2) \to \mathbb{R}$, we denote $f (Z^i_t, \Pi^N_t)$ and $f (\tilde Z^i_t, \tilde \Pi^N_t)$ by $f^i_t$ and $\tilde f^i_t$ respectively. Let us detail only the partial observation case. As in the proof of Theorem \ref{th: convergence InN} we start by detailing the convergence of $\tilde{I}_n^{N,P}(f)$, defined as in \eqref{eq: def I tilde}, and then we deduce the result for ${I}_n^{N,P}(f)$. 
Using the decomposition $\tilde X^i_{(j+1)\Delta}- \tilde X^i_{j\Delta} - \Delta \tilde b^i_{j\Delta} = \Delta^{\frac{1}{2}} a^i_{j\Delta} (\tilde \xi_j + \xi^i_{j+1}) + \tilde \varepsilon^i_j$ as in Proposition \ref{prop: increments tilde}(i), we rewrite $\tilde I^{N,P}_n (f_\theta)$ and reorder its terms as follows:
\begin{align*}
N^{-\frac{1}{2}} \tilde{I}_n^{N,P}(f) =
F_0 + \sum_{j = 2}^{n -2} F_{j,1} + \sum_{j = 1}^{n -2}  F_{j,2},
\end{align*}
where 
\begin{align*}
F_0 &:= \Big( \frac{\Delta}{N} \Big)^{\frac 1 2} \sum_{i=1}^N ( \tilde f_0^i a_{\Delta}^i \tilde \xi^i_1 + \tilde f_{(n-3)\Delta}^i a_{(n-2)\Delta}^i \xi^i_{n-1}),\\
F_{j,1} &:= \Big( \frac{\Delta}{N} \Big)^{\frac 1 2} \sum_{i=1}^N ( \tilde f_{(j-2)\Delta}^i a_{(j-1)\Delta}^i \xi^i_j + \tilde f_{(j-1)\Delta}^i a_{j\Delta}^i \tilde \xi^i_j),\\
F_{j,2} &:= N^{-\frac{1}{2}} \sum_{i=1}^N \tilde f_{(j-1)\Delta}^i \tilde \varepsilon^i_j.
\end{align*}
We want to prove that 
\begin{equation}\label{lim:ANF}
\sum_{j=2}^{n-2} F_{j,1} \xrightarrow{d} \mathcal{N}(0, \bar{\Pi} ( (fa)^2 ) ), \qquad \sum_{j=1}^{n-2} F_{j,2} \xrightarrow{\P} 0, \qquad F_0 \xrightarrow{\P} 0.
\end{equation}

Our proof of the first relation uses a central limit theorem for martingale difference arrays (see, e.g.\ Theorems 3.2 and 3.4 of \cite{38McK}). In particular, we need to prove that  
\begin{align}
\sum_{j = 2}^{n -2} \E_{j\Delta}[F_{j,1}] &\xrightarrow{\mathbb{P}} 0,\label{eq: conv Zbar n+9}\\
\sum_{j = 2}^{n -2} \E_{j\Delta}[(F_{j,1})^2] &\xrightarrow{\mathbb{P}} \bar{\Pi} ( (f a)^2 ),\label{eq: conv Zbar2 n+10}\\
\sum_{j = 2}^{n -2} \E_{j\Delta}[(F_{j,1})^4] &\xrightarrow{\mathbb{P}}0.\label{eq: conv Zbar2+r n+11}
\end{align}
The validity of \eqref{eq: conv Zbar n+9} is straightforward from $\E_{j\Delta} [\xi^i_j] = \E_{j\Delta} [\tilde \xi_j^i] = 0$. Regarding \eqref{eq: conv Zbar2 n+10}, we have that $\E_{j\Delta}[\xi^{i_1}_j \xi^{i_2}_j] = \E_{j\Delta}[\tilde \xi^{i_1}_j \tilde \xi^{i_2}_j] = 2\E_{j\Delta}[\xi^{i_1}_j \tilde \xi^{i_2}_j] = \frac{1}{3}$ if $i_1 = i_2$, whereas all these terms are zero if $i_1 \neq i_2$, providing
\begin{align*}
\sum_{j=2}^{n-2} \E_{j\Delta}[(F_{j,1})^2] &= \frac{\Delta}{3 N} \sum_{j=2}^{n-2} \sum_{i = 1}^N (( \tilde f_{(j-2)\Delta}^i a_{(j-1)\Delta}^i)^2\\
&\qquad \qquad + \tilde f_{(j-2)\Delta}^i a_{(j-1)\Delta}^i \tilde f_{(j-1)\Delta}^i a_{j\Delta}^i + (\tilde f_{(j-1)\Delta}^i a_{j\Delta}^i )^2).
\end{align*}
Its convergence to $\bar{\Pi} ( (f a)^2 )$ follows in the same ways as the first convergence in \eqref{lim:omega1}. Let us turn to \eqref{eq: conv Zbar2+r n+11}. WLOG assume that $1 \le i_1 \le \dots \le i_4 \le N$ and $w_k \in \{0,1\}$, $k=1,\dots,4$. Then 
$$
\E_{j\Delta} \Big[ \prod_{k=1}^4 (w_k \xi^{i_k}_j + (1-w_k) \tilde \xi^{i_k}_j) \Big] = C,
$$ 
where $C \neq 0$ only if $i_1=i_2$ and $i_3=i_4$. Hence, the $L^1$-norm of the LHS of \eqref{eq: conv Zbar2+r n+11} is bounded by $C \Delta^{-1} (N^2 +N) \Delta^2 N^{-2} \le C \Delta$, which tends to zero, proving \eqref{eq: conv Zbar2+r n+11}.


The second relation in \eqref{lim:ANF} follows by \cite[Lemma 9]{GenJac93} if we prove that both
\begin{align*}
\Sigma_1 := \sum_{j=1}^{n-2} \E_{j\Delta} [F_{j,2}] 
&= N^{-\frac{1}{2}} \sum_{j=1}^{n-2} \sum_{i=1}^N \tilde f^i_{(j-1)\Delta} \E_{j\Delta} [\tilde \varepsilon^i_j],\\ 
\Sigma_2:= \sum_{j=1}^{n-2} \E_{j\Delta} [(F_{j,2})^2] 
&= N^{-1} \sum_{j=1}^{n-2} \sum_{i_1,i_2=1}^N  \tilde f^{i_1}_{(j-1)\Delta} \tilde f^{i_2}_{(j-1)\Delta} \E_{j\Delta} [\tilde \varepsilon^{i_1}_j \tilde \varepsilon^{i_2}_j]
\end{align*}
converge in $L^1$ to $0$. We note that $\E [|\E_{j\Delta}[\tilde \varepsilon^i_j]|^k]^{\frac{1}{k}} \le C \Delta^{\frac{3}{2}}$ and $\E [|\tilde \varepsilon^i_j|^k]^{\frac{1}{k}} \le C \Delta$ for all $k\ge1$ by Proposition \ref{prop: increments tilde}(i). Hence, we have that $\E[|\Sigma_1|] \le C (N\Delta)^{\frac{1}{2}}$ and $\E[|\Sigma_2|] \le C N\Delta$, which tend to $0$. Also, the last relation in \eqref{lim:ANF} holds true, since $\E [|F_0|] \le C \Delta^{\frac{1}{2}} N^{\frac{1}{2}}$. 


This gives the convergence in distribution of $N^{-\frac{1}{2}} \tilde{I}_n^{N,P}(f)$. To deduce the result for $N^{-\frac{1}{2}} {I}_n^{N,P}(f)$, it is enough to show that
\begin{equation}{\label{eq: tilde to not INn n + 12}}  
N^{-\frac{1}{2}} (I^{N,P}_n (f)-\tilde I^{N,P}_n(f)) = \frac{\Delta}{N^{\frac{1}{2}}} \sum_{j = 1}^{n -2} \sum_{i = 1}^N \tilde f_{(j- 1)\Delta}^i (\tilde b_{j\Delta}^i - \tilde b_{(j-1)\Delta}^i) \xrightarrow{L^1} 0.
\end{equation}
We recall that $f$ has a polynomial growth, whereas $b$ is Lipschitz continuous according to \ref{as2}. Using the Cauchy-Schwarz inequality and Proposition \ref{prop: increments tilde}(ii), the $L^1$-norm of the RHS of \eqref{eq: tilde to not INn n + 12} gets bounded by $C (N\Delta)^{\frac 1 2}$.
It tends to $0$ under our hypothesis, that concludes the proof of our result.

\end{proof}

\subsection{Proof of Theorem \ref{th: convergence law QnN}}

\begin{proof}
Similarly as before, we only provide the details for the partial observation case. Remark that one should closely follow the proof of (38) in \cite{McKean} to see the variance term appearing, in the complete observation case. In the partial observation framework we have that 
$$
(N \Delta)^{-\frac{1}{2}} (Q^{N,P}_n (f) - \frac{2}{3} \Delta \nu_n^{N,P} (f a^2)) = \sum_{k=1}^5 R^N_{n,k},
$$
where by the definitions of $Q^{N,P}_n$, $\nu^{N,P}_n$ and Proposition \ref{prop: increments tilde}(i),
\begin{align*}
R^N_{n,1} &= \Big( \frac{\Delta}{N} \Big)^{\frac{1}{2}} \sum_{j=1}^{n-2} \sum_{i=1}^N (f a^2)^i_{j\Delta} \Big( (U^i_j)^2 - \frac{2}{3} \Big)\\
R^N_{n,2} &= \Big( \frac{\Delta}{N} \Big)^{\frac{1}{2}} \sum_{j=1}^{n-2} \sum_{i=1}^N ( {\tilde f}^i_{(j-1)\Delta} - f^i_{j\Delta}) (a^2)^i_{j\Delta} (U^i_j)^2\\
R^N_{n,3}&=  \frac{2}{3} \Big( \frac{\Delta}{N} \Big)^{\frac{1}{2}} \sum_{j=1}^{n-2} \sum_{i=1}^N ( f^i_{j\Delta} - \tilde f^i_{j\Delta} ) (a^2)^i_{j\Delta} ,\\
R^N_{n,4} &= - \frac{2}{3} \Big( \frac{\Delta}{N} \Big)^{\frac{1}{2}}\sum_{i=1}^N (\tilde f^i_0  (a^2)^i_0 + \tilde f^i_{(n-1)\Delta} (a^2)^i_{(n-1)\Delta}), \\
R^N_{n,5} &= (N\Delta)^{-\frac{1}{2}} \sum_{j=1}^{n-2} \sum_{i=1}^N \tilde f^i_{(j-1)\Delta} (2\Delta^{\frac{1}{2}} a^i_{j\Delta} U^i_j + \Delta \tilde b^i_{j\Delta} + \tilde \varepsilon^i_j ) (\Delta \tilde b^i_{j\Delta} + \tilde \varepsilon^i_j).
\end{align*}
We will prove that as $N\Delta \to 0$,
\begin{equation}\label{lim:RNn}
R^N_{n,1} \xrightarrow{d} \mathcal{N}(0,\bar \Pi(f^2a^4)), \qquad \sum_{k=2}^5 R^N_{n,k} \xrightarrow{\P} 0.
\end{equation}
We first study the convergence of $R^N_{n,1}$. We recall that $U^i_j = \tilde \xi^i_j + \xi^i_{j+1}$, where $\tilde \xi^i_j$, $\xi^i_{j+1}$ are $\mathcal{F}^N_{(j+1)\Delta}$-, $\mathcal{F}^N_{(j+2)\Delta}$-measurable, respectively, and independent, and $\mathcal{N}(0,\frac{1}{3})$-distributed. We rewrite
\begin{equation*}
R^N_{n,1} = \sum_{j=0}^{n-2} \rho_{j},
\end{equation*}
in order to obtain a martingale difference array 
for $j=1,\dots,n-2$:
\begin{align*}
\rho_j &:= \Big( \frac{\Delta}{N} \Big)^{\frac{1}{2}} \sum_{i=1}^N \Big( (f a^2)^i_{(j-1)\Delta} \Big( (\xi^i_j)^2-\frac{1}{3} \Big) + (f a^2)^i_{j\Delta} \Big( (\tilde \xi^i_j)^2 - \frac{1}{3} \Big) + 2 (f a^2)^i_{(j-1)\Delta} \tilde \xi^i_{j-1} \xi^i_j \Big),\\
\rho_0 &:= \Big( \frac{\Delta}{N} \Big)^{\frac{1}{2}} \sum_{i=1}^N \Big( (f a^2)^i_\Delta \Big( (\tilde \xi^i_0)^2-\frac{1}{3} \Big) + (f a^2)^i_{(n-2)\Delta} \Big( (\xi^i_{n-1})^2-\frac{1}{3} + 2\tilde \xi^i_{n-2} \xi^i_{n-1} \Big) \Big).
\end{align*}
We get that $\E[|\rho_0|] \le C (N\Delta)^{\frac{1}{2}}$, which tends to $0$, while the sum of $\rho_j$ over $j=1,\dots,n-2$ is the main term. 
In order to prove that this sum converges in distribution to $\mathcal{N}(0,\bar \Pi(f^2a^4))$, we use a central limit theorem for martingale difference arrays (see, e.g.\ Theorems 3.2 and 3.4 of \cite{38McK}). In particular, we need to prove that 
\begin{align}{\label{eq: conv omega bar n+13}}
\sum_{j = 1}^{n -2} \E_{j\Delta}[\rho_j] &\xrightarrow{\mathbb{P}}0, \\
{\label{eq: conv omega bar squared n+14}}
\sum_{j = 1}^{n -2} \E_{j\Delta}[(\rho_j)^2] &\xrightarrow{\mathbb{P}} \bar{\Pi}(f^2 a^4),\\
{\label{eq: conv omega bar r n+15}}
\sum_{j = 1}^{n -2} \E_{j\Delta}[(\rho_j)^4] &\xrightarrow{\mathbb{P}} 0. 
\end{align}
The validity of \eqref{eq: conv omega bar n+13} is straightforward from $\E_{j\Delta}[\rho_{j}]  = 0$. 
Next, we calculate that 
\begin{align*}
\E_{j\Delta}[(\rho_{j})^2]
&= \frac{\Delta}{N} \sum_{i=1}^N \Big( \frac{2}{9} (f^2 a^4)^i_{(j-1)\Delta} + \frac{2}{9} (f^2 a^4)^i_{j\Delta}\\ 
&\qquad + \frac{4}{3} (f^2a^4)^i_{(j-1)\Delta} (\tilde \xi^i_{j-1})^2 + \frac{1}{9} (f a^2)^i_{(j-1)\Delta}(f a^2)^i_{j\Delta} \Big),
\end{align*}
where 
$$
\frac{\Delta}{N} \sum_{j=1}^{n-2} \sum_{i=1}^N (f^2 a^4)^i_{(j-2)\Delta} (\tilde \xi^i_{j-1})^2 \xrightarrow{\P} \frac{1}{3} \bar \Pi (f^2 a^4)
$$
follows by \cite[Lemma 9]{GenJac93} (see also \eqref{lim:omega1}) and by Proposition \ref{prop: converence nu}. This concludes the proof of \eqref{eq: conv omega bar squared n+14}. Finally, \eqref{eq: conv omega bar r n+15} follows from the bound $C \Delta \to 0$ on the $L^1$-norm of its LHS (see also \eqref{eq: conv Zbar2+r n+11}).

It remains to prove the relation on the RHS of \eqref{lim:RNn}. Let us first prove that $R^N_{n,2} \xrightarrow{\P} 0$ as $N\Delta \to 0$. We note that
$$
\Big( \frac{\Delta}{N} \Big)^{\frac 1 2} \sum_{j=1}^{n-2} \sum_{i=1}^N (\tilde f_{(j-1)\Delta} - f^i_{j\Delta})((a^2)^i_{j\Delta} - (a^2)^i_{(j-1)\Delta}) (U^i_j)^2 \xrightarrow{L^1} 0
$$
as $N\Delta \to 0$, because the condition \eqref{eq: cond pol growth} on $f, a^2$ and Lemma \ref{l: moments}, Proposition \ref{p: error X tilde}(ii), (iii) give us
$$
\E [| (\tilde f_{(j-1)\Delta} - f^i_{j\Delta})((a^2)^i_{j\Delta} - (a^2)^i_{(j-1)\Delta}) (U^i_j)^2 |] \le C \Delta.
$$
Furthermore, application of \cite[Lemma 9]{GenJac93} implies that 
$$
\Big( \frac{\Delta}{N} \Big)^{\frac{1}{2}} \sum_{j=1}^{n-2} \sum_{i=1}^N (\tilde f_{(j-1)\Delta} - f^i_{j\Delta}) (a^2)^i_{(j-1)\Delta} (U^i_j)^2 \xrightarrow{\P} 0
$$
as $N \Delta \to 0$, because using Proposition \ref{prop: ito} and Parts (ii), (iii) of Lemma \ref{l: moments}, Proposition \ref{p: error X tilde}
we can bound the $L^1$-norm of 
$$
\E_{(j-1)\Delta}[(\tilde f^i_{(j-1)\Delta} - f^i_{j\Delta}) (a^2)^i_{(j-1)\Delta} (U^i_j)^2] = \frac{2}{3} (a^2)^i_{(j-1)\Delta} \E_{(j-1)\Delta}[\tilde f^i_{(j-1)\Delta} - f^i_{j\Delta}]
$$
and
\begin{align*}
    &\E_{(j-1)\Delta}[(\tilde f^{i_1}_{(j-1)\Delta} - f^{i_1}_{j\Delta}) (a^2)^{i_1}_{(j-1)\Delta} (U^{i_1}_j)^2(\tilde f^{i_2}_{(j-1)\Delta} - f^{i_2}_{j\Delta}) (a^2)^{i_2}_{(j-1)\Delta} (U^{i_2}_j)^2] \\
    &\qquad = C (a^2)^{i_1}_{(j-1)\Delta} (a^2)^{i_2}_{(j-1)\Delta} \E_{(j-1)\Delta}[(\tilde f^{i_1}_{(j-1)\Delta} - f^{i_1}_{j\Delta})(\tilde f^{i_2}_{(j-1)\Delta} - f^{i_2}_{j\Delta})]
\end{align*}
by $C \Delta$ and $C \Delta$, respectively. Similarly, application of \cite[Lemma 9]{GenJac93} implies that $R^N_{n,3} \xrightarrow{\P} 0$ as $N\Delta \to 0$. Also, it is easy to check that $\E [|R^N_{n,4}|] \le C (N \Delta)^{\frac{1}{2}}$ tends to $0$. At last, let us prove that $R^N_{n,5} \xrightarrow{\P} 0$ as $N\Delta \to 0$. 
By applying the $L^1$ argument again, the proof reduces to that of
$$
N^{-\frac{1}{2}} \sum_{j=1}^{n-2} \sum_{i=1}^N (f a^2)^i_{j\Delta} U^i_j (\Delta b^i_{j\Delta} + \tilde \varepsilon^i_j) \xrightarrow{\P} 0
$$
as $N\Delta \to 0$, which in turn follows by \cite[Lemma 9]{GenJac93}, because using Proposition \ref{prop: increments tilde}(i) we can bound the $L^1$-norm of
$$
\E_{j\Delta} [(f a^2)^i_{j\Delta} U^i_j (\Delta b^i_{j\Delta} + \tilde \varepsilon^i_j) ] = (f a^2)^i_{j\Delta} \E_{j\Delta} [U^i_j \tilde \varepsilon^i_j]
$$
and  
$$
\E_{j\Delta} [(f a^2)^{i_1}_{j\Delta} U^{i_1}_j (\Delta b^{i_1}_{j\Delta} + \tilde \varepsilon^{i_1}_j) (f a^2)^{i_2}_{j\Delta} U^{i_2}_j (\Delta b^{i_2}_{j\Delta} + \tilde \varepsilon^{i_2}_j) ]
$$
by $C \Delta^{\frac{3}{2}}$ and $C\Delta^2$, respectively. This completes the proof of the relation on the RHS of \eqref{lim:RNn} and that of the theorem.
\end{proof}

\subsection{Proof of Proposition \ref{prop: hypo}}
\begin{proof}
In the proof below, we follow the formalism introduced in \cite[Section 2.3.2]{Nualart}. Note that since the algebraic vector corresponding to the vector field $A_{2k}$ is given by ${\spa \hat{A}}_{2k} = (0, \dots, a^{(k)}(Z_t), \dots, 0)^\top$, and $a^{(k)}(Z_t)$ is bounded away from zero, the vectors $({\spa \hat{A}}_{2k})_{k = 1, \dots, N}$ are linearly independent and span an $N$-dimensional subspace of $\mathbb{R}^{2N}$. Our goal is to show that to generate the full tangent space $\mathbb{R}^{2N}$, it is sufficient to include the coordinate vectors of the Lie brackets, {\spa denoted by $\widehat{[A_0, A_{2k}]}$, for $k = 1, \dots, N$.}

In our framework, the sparsity of the diffusion matrix $\hat{A}$ significantly simplifies the analysis. Recall the definition of the drift vector field:
$$A_0 = B - \frac{1}{2} \sum_{l=1}^{2N} A_l^\nabla A_l.$$
Observe that this reduces to:
$$A_0 = B - \frac{1}{2} \sum_{l=1}^{2N} A_{2l}^\nabla A_{2l},$$
since $A_l^\nabla A_l = 0$ for any odd index $l$ (due to ${\spa \hat{A}}_{2l-1} = \mathbf{0}$). Thus, only the even-indexed terms contribute. Moreover, according to the definition of the covariant derivative in \eqref{eq: cov der}, we have:
$$A_{2l}^\nabla A_{2l} = \sum_{i,j = 1}^{2N} {\spa \hat{A}}^i_{2l} \partial_{z_i} {\spa \hat{A}}_{2l}^j \partial_{z_j},$$
where $\partial_{z_j} = \mathbf{1}_{\{j = 2k -1\}} \partial_{y^k} + \mathbf{1}_{\{j = 2k\}} \partial_{x^k}$. Recall also that for all $i \neq 2l$, we have ${\spa \hat{A}}^i_{2l} = 0$, while ${\spa \hat{A}}_{2l}^{2l}= a^{(l)}(Z_t)$. Therefore, the operator simplifies to:
$$A_{2l}^\nabla A_{2l} = a^{(l)}(Z_t) \partial_{x^l} a^{(l)}(Z_t) \partial_{z_{2l}}.$$

Let ${\spa \hat{A}}_0$ denote the algebraic column vector containing the coefficients of the vector field $A_0$. From the derivation above, we can write ${\spa \hat{A}}_0$ explicitly as:
\begin{equation*}
{\spa \hat{A}}_0 = 
\begin{pmatrix}
 b_1^{(1)}(Z_t) \\
 b_2^{(1)}(Z_t) \\
 \vdots \\
 b_1^{(N)}(Z_t) \\
 b_2^{(N)}(Z_t)
\end{pmatrix} 
- \frac{1}{2}
\begin{pmatrix}
 0 \\
 a^{(1)}(Z_t) \partial_{x^1} a^{(1)}(Z_t) \\
 \vdots \\
 0 \\
 a^{(N)}(Z_t) \partial_{x^N} a^{(N)}(Z_t)
\end{pmatrix}
=
\begin{pmatrix}
 b_1^{(1)}(Z_t) \\
 b_2^{(1)}(Z_t) - \frac{1}{2} a^{(1)}(Z_t) \partial_{x^1} a^{(1)}(Z_t) \\
 \vdots \\
 b_1^{(N)}(Z_t) \\
 b_2^{(N)}(Z_t) - \frac{1}{2} a^{(N)}(Z_t) \partial_{x^N} a^{(N)}(Z_t)
\end{pmatrix}.
\end{equation*}

We can now compute the Lie brackets $[A_0, A_{2k}]$ for $k = 1, \dots, N$. By definition:
$$[A_0, A_{2k}] = A_0^\nabla A_{2k} - A_{2k}^\nabla A_0.$$
For the first term, we have:
$$A_0^\nabla A_{2k} = \sum_{i= 1}^{2N} {\spa \hat{A}}_0^i \partial_{z_i} a^{(k)}(Z_t) \partial_{z_{2k}},$$
which corresponds to the algebraic vector:
$$\widehat{A_0^\nabla A_{2k}} = 
\begin{pmatrix}
 0 \\
 \vdots \\
 \sum_{i= 1}^{2N} {\spa \hat{A}}_0^i \partial_{z_i} a^{(k)}(Z_t) \\
 \vdots \\
 0
\end{pmatrix},$$
where the only non-zero component is in the $2k$-th position. For the second term, we have:
$$A_{2k}^\nabla A_0 = \sum_{j = 1}^{2N} a^{(k)}(Z_t) \partial_{x^{k}} {\spa \hat{A}}_0^j \partial_{z_j}.$$
Let us introduce, for $l = 1, \dots , N$, the function $f_l(Z_t) := b_2^{(l)}(Z_t) - \frac{1}{2} a^{(l)}(Z_t) \partial_{x^l} a^{(l)}(Z_t)$. The coordinate vector for the Lie bracket, {\spa $\widehat{[A_0, A_{2k}]}$}, is then given by:
$$ {\spa \widehat{[A_0, A_{2k}]}} =
\begin{pmatrix}
- a^{(k)}(Z_t) \partial_{x^{k}} {\spa \hat{A}}_0^1 \\
- a^{(k)}(Z_t) \partial_{x^{k}} {\spa \hat{A}}_0^2 \\
 \vdots \\
 - a^{(k)}(Z_t) \partial_{x^{k}} {\spa \hat{A}}_0^{2k - 1} \\
- a^{(k)}(Z_t) \partial_{x^{k}} {\spa \hat{A}}_0^{2k} + \sum_{i= 1}^{2N} {\spa \hat{A}}_0^i \partial_{z_i} a^{(k)}(Z_t)  \\ 
\vdots \\
- a^{(k)}(Z_t) \partial_{x^{k}} {\spa \hat{A}}_0^{2N - 1} \\
- a^{(k)}(Z_t) \partial_{x^{k}} {\spa \hat{A}}_0^{2N}
\end{pmatrix} 
= 
\begin{pmatrix}
- a^{(k)}(Z_t) \partial_{x^{k}} b_1^{(1)}(Z_t) \\
- a^{(k)}(Z_t) \partial_{x^{k}} f_1(Z_t) \\
 \vdots \\
 - a^{(k)}(Z_t) \partial_{x^{k}} b_1^{(k)}(Z_t) \\
- a^{(k)}(Z_t) \partial_{x^{k}} f_k(Z_t) + \sum_{i= 1}^{2N} {\spa \hat{A}}_0^i \partial_{z_i} a^{(k)}(Z_t)  \\ 
\vdots \\
- a^{(k)}(Z_t) \partial_{x^{k}} b_1^{(N)}(Z_t) \\
- a^{(k)}(Z_t) \partial_{x^{k}} f_N(Z_t)
\end{pmatrix}.
$$

Observe that the vectors ${\spa \widehat{[A_0, A_{2k_1}]}}$ and ${\spa \widehat{[A_0, A_{2k_2}]}}$ are linearly independent for $k_1 \neq k_2$. This follows from the assumption that for any $j, k \in \{1, \dots, N \}$ with $j \neq k$, we have $\partial_{x^k} b_1^{(k)}(z) \neq \partial_{x^k} b_1^{(j)}(z) \neq 0$. 

Consequently, the basis for $\mathbb{R}^{2N}$ is given by the algebraic vectors ${\spa (\hat{A}_{2k}, \widehat{[A_0, A_{2k}]})_{k \in \{1, \dots, N \}}}$. This holds precisely because we assumed $\partial_{x^{k}} b_1^{(i)}(z) \neq 0$ for all $k, i \in \{1, \dots, N \}$, and that the diffusion coefficient $a^{(k)}(z)$ is non-degenerate for each $k \in \{1, \dots, N \}$. Thus, Hörmander's condition is satisfied, and the proof is complete.

\end{proof}

\end{document}